\documentclass[reqno]{amsart}
\usepackage{amsmath}%
\usepackage{amsfonts}%
\usepackage{amssymb}%
\usepackage{graphicx}%
\usepackage[latin9]{inputenc}
\usepackage[T1]{fontenc}     
\usepackage{geometry}         
\usepackage[dvips]{color}
\usepackage{amsmath, latexsym, amsfonts, amssymb, amsthm, amscd}
\usepackage{graphics,epsf,psfrag,epsfig}
\numberwithin{equation}{section}

\newtheorem{theorem}{Theorem}[section]
\newtheorem{lemma}[theorem]{Lemma}
\newtheorem{proposition}[theorem]{Proposition}

\newtheorem{rem}[theorem]{Remark}
\newtheorem{definition}[theorem]{Definition}
\newtheorem{hyp}[theorem]{Assumption}

\newcommand{\ind}{\mathbf{1}}

\newcommand{\R}{\mathbb{R}}
\newcommand{\Z}{\mathbb{Z}}
\newcommand{\N}{\mathbb{N}}
\renewcommand{\tilde}{\widetilde}
\renewcommand{\hat}{\widehat}

\DeclareMathSymbol{\leqslant}{\mathalpha}{AMSa}{"36} 
\DeclareMathSymbol{\geqslant}{\mathalpha}{AMSa}{"3E} 
\DeclareMathSymbol{\eset}{\mathalpha}{AMSb}{"3F}     
\renewcommand{\leq}{\;\leqslant\;}                   
\newcommand{\dd}{\,\text{\rm d}}             

\newcommand{\sumtwo}[2]{\sum_{\substack{#1 \\ #2}}} 

\newcommand{\cN}{{\ensuremath{\mathcal N}} }

\newcommand{\cZ}{{\ensuremath{\mathcal Z}} }
\newcommand{\cI}{{\ensuremath{\mathcal I}} }

\newcommand{\bP}{{\ensuremath{\mathbf P}} }
\newcommand{\bE}{{\ensuremath{\mathbf E}} }

\newcommand{\bbE}{{\ensuremath{\mathbb E}} }

\newcommand{\bbP}{{\ensuremath{\mathbb P}} }

\newcommand{\ga}{\alpha}
\newcommand{\gb}{\beta}
\newcommand{\gga}{\gamma}            
\newcommand{\gd}{\delta}
\newcommand{\gep}{\varepsilon}       

\newcommand{\gz}{\zeta}

\newcommand{\gD}{\Delta}
\newcommand{\gk}{\kappa}

\newcommand{\gl}{\lambda}

\newcommand{\gs}{\sigma}
\newcommand{\gS}{\Sigma}

\newcommand{\gba}{\gb_c^{ann}}
\newcommand{\gbc}{\gb_c}

\renewcommand{\P}{{\ensuremath{\mathbf P}} }
\newcommand{\E}{{\ensuremath{\mathbf E}} }
\newcommand{\PI}{P_{\cI}}

\newcommand{\Znb}{Z_{N,Y}^{\gb}}

\newcommand{\Py}{{\mathbb P}^Y}

\newcommand{\Ey}{{\mathbb E}^Y}
\newcommand{\Px}{{\mathbb P}^X}
\newcommand{\px}{p^X}
\newcommand{\Pnb}{{\mathbb P}_{N,Y}^{\gb}}
\newcommand{\Enb}{{\mathbb E}_{N,Y}^{\gb}}

\newcommand{\Ex}{\mathbb E^X}
\newcommand{\Pxy}{{\mathbb P}^{X-Y}}
\newcommand{\Exy}{{\mathbb E}^{X-Y}}

\newcommand{\norm}[1]{\left\| #1 \right\|}

\newcommand{\Ninf}{N\to\infty}

\newcommand{\Fa}{F^{ann}}
\newcommand{\xyn}{X_n=Y_n}
\newcommand{\xyN}{X_N=Y_N}
\newcommand{\ann}{\textit{annealed}}

\newcommand{\Znzc}{\check Z_{N,Y}^{z}}
\newcommand{\Gxy}{G^{X-Y}}
\newcommand{\pxy}{p^{X-Y}}

\newcommand{\ZI}{Z_{z,Y}^{\cI}}

\newcommand{\EL}{{\ensuremath{\hat \bE}}}

\newcommand{\PL}{{\ensuremath{\hat \bP}}}
\newcommand{\Pct}{{\ensuremath{\mathbb P}}_{\tau}}
\newcommand{\Ectdf}{\mathbb{E}_{\tau}}
\newcommand{\Ect}{\mathbb{E}_{\tau}}
\newcommand{\gdt}{\gD\tau}
\newcommand{\EI}{E_{\cI} }
\newcommand{\Eyrx}{\mathbb{E}_{r,x}^{Y}}
\newcommand{\Pyrx}{\mathbb{P}_{r,x}^{Y}}
\newcommand{\Edf}{\hat\bE}
\newcommand{\Pdf}{\hat\bP}

\begin{document}

\title[Random Walk Pinning model in $d=3$]{On the critical point of the \\
Random Walk Pinning Model in dimension $d=3$}

\author{Quentin Berger}
\address{
Laboratoire de Physique, ENS Lyon,  Universit\'e de Lyon, 46 All\'ee d'Italie, 
69364 Lyon, France
}
\email{quentin.berger@ens.fr}
\author{Fabio Lucio Toninelli}
\address{CNRS and 
Laboratoire de Physique, ENS Lyon, Universit\'e de Lyon, 46 All\'ee d'Italie, 
69364 Lyon, France
}
\email{fabio-lucio.toninelli@ens-lyon.fr}

 \thanks{This work was
  supported by the European Research Council through the ``Advanced
  Grant'' PTRELSS 228032, and by ANR through the grant LHMSHE}

\begin{abstract}
  We consider the Random Walk Pinning Model studied in \cite{BS08} and
  \cite{BGdH08}: this is a random walk $X$ on $\Z^d$, whose law is
  modified by the exponential of $\beta$ times $L_N(X,Y)$, the
  collision local time up to time $N$ with the (quenched) trajectory
  $Y$ of another $d$-dimensional random walk. If $\beta$ exceeds a
  certain critical value $\beta_c$, the two walks stick together for
  typical $Y$ realizations (localized phase).
A natural question is whether the disorder is relevant or not, that is whether the
  {\sl quenched} and {\sl annealed} systems have the same critical behavior.
  Birkner and Sun
  \cite{BS08} proved that $\beta_c$ coincides with the critical point
  of the {\sl annealed} Random Walk Pinning Model if the space dimension is
  $d=1$ or $d=2$, and that it differs from it in dimension $d\ge4$
  (for $d\ge 5$, the result was proven also in \cite{BGdH08}). Here,
  we consider the open case of the {\sl marginal dimension} $d=3$, and
  we prove non-coincidence of the critical points.
\end{abstract}

\maketitle

\section{Introduction}
We consider the Random Walk Pinning Model (RWPM): the starting point is a
zero-drift random walk $X$ on $\Z^d$ ($d\ge1$), whose law is modified
by the presence of a second random walk, $Y$. The trajectory of $Y$ is
fixed (quenched disorder) and can be seen as the random medium. The
modification of the law of $X$ due to the presence of $Y$ takes the
Boltzmann-Gibbs form of the exponential of a certain interaction
parameter, $\beta$, times the collision local time of $X$ and $Y$
up to time $N$, $L_N(X,Y):=\sum_{1\le n\le N}\ind_{\{X_n=Y_n\}}$.
If $\beta$ exceeds a certain threshold value $\beta_c$, then for almost every
realization of $Y$ the walk $X$ sticks together with $Y$, in the 
thermodynamic limit $N\to\infty$. If on the other hand $\beta<\beta_c$, 
then $L_N(X,Y)$ is $o(N)$ for typical trajectories.

Averaging with respect to $Y$ the partition function, one obtains the partition function of the
so-called annealed model, whose critical point $\beta_c^{ann}$ is
easily computed; a natural question is whether
$\beta_c\ne\beta_c^{ann}$ or not.  In the renormalization group
language, this is related to the question whether disorder is {\sl
  relevant} or not.  In an early version of the paper \cite{BGdH08},
Birkner {\sl et al.}  proved that $\beta_c\ne\beta_c^{ann}$ in
dimension $d\ge 5$. Around the same time, Birkner and Sun \cite{BS08}
extended this result to $d=4$, and also proved that the two critical
points {\sl do coincide} in dimensions $d=1$ and $d=2$. 

The dimension $d=3$ is the {\sl marginal dimension} in the
renormalization group sense, where not even heuristic arguments like the
``Harris criterion'' (at least its most naive version) can predict whether one has disorder relevance
or irrelevance.  Our main result here is that quenched and annealed
critical points differ also in $d=3$. 

For a discussion of the connection of the RWPM with the ``parabolic
Anderson model with a single catalyst'', and of the implications of
$\beta_c\ne\beta_c^{ann}$ about  the location of the weak-to-strong transition for the directed polymer in
random environment, we refer to \cite[Sec. 1.2 and 1.4]{BS08}.

\medskip

Our proof is based on the idea of bounding the fractional moments of
the partition function, together with a suitable change of measure
argument.  This technique, originally introduced in
\cite{DGLT07,GLT08,GLT09} for the proof of disorder relevance for the
random pinning model with tail exponent $\ga\ge 1/2$, has also proven
to be quite powerful in other cases: in the proof of
non-coincidence of critical points for the RWPM in dimension $d\ge 4$
\cite{BS08}, in the proof that ``disorder is always strong'' for the
directed polymer in random environment in dimension $(1+2)$
\cite{Lacoin} and finally in the proof that quenched and annealed
large deviation functionals for random walks in random environments in
two and three dimensions differ \cite{YZ}.
Let us mention that for the random pinning model there is another method, developed by Alexander and Zygouras
\cite{AZ08}, to prove disorder relevance: however, their method fails in the marginal situation $\alpha=1/2$ (which 
corresponds to $d=3$ for the RWPM).

\medskip

To guide the reader through the paper, let us point out immediately
what are the  novelties and the similarities of our proof with respect
to the previous applications of the fractional moment/change of measure 
method: 

\begin{itemize}

\item the change of measure chosen by Birkner and Sun in \cite{BS08}
consists essentially in correlating positively each increment
of the random walk $Y$ with the next one. Therefore, under the modified measure,
$Y$ is more diffusive. The change of measure we use in dimension three 
has also  the effect of correlating positively the increments of $Y$, but in our case
the correlations have long range (the correlation between the $i^{th}$ and the $j^{th}$ increment decays like $|i-j|^{-1/2}$).
Another ingredient which was absent in \cite{BS08} and which is essential
in $d=3$ is a coarse-graining step, of the type of that employed 
in \cite{T_cg,GLT09};

\item while the scheme of the proof of our Theorem \ref{th:main}
has many points in common with  that of \cite[Th. 1.7]{GLT09}, here
we need new renewal-type estimates  (e.g. Lemma \ref{th:derniert})
and a careful application of the Local Limit Theorem to prove that the
average of the partition function under the modified measure is small
(Lemmas \ref{EctF big} and \ref{variance}).

\end{itemize}

\section{Model and results}

\subsection{The random walk pinning model}

Let $X=\{X_n\}_{n\geq 0}$ and $Y=\{Y_n\}_{n\geq 0}$ be two independent
discrete-time random walks on $\Z^d$, $d\ge 1$, starting from $0$, and
let $\Px$ and $\Py$ denote their respective laws.  We make the
following assumption:
\begin{hyp}
\label{hypwright}
\rm The random walk $X$ is aperiodic.
The increments $(X_i-X_{i-1})_{i\ge1}$ are i.i.d., symmetric and 
have a sub-Gaussian tail:  for every $x>0$,
\begin{eqnarray}
  \label{eq:hypwright}
  \bbP^X(\|X_1\|\ge x)\le M\int_x^\infty \exp(-h t^2)\dd t
\end{eqnarray}
for some positive constants $M,h$, where $\|\cdot\|$ denotes the
Euclidean norm on $\Z^d$. Moreover, the covariance matrix of $X_1$,
call it $\Sigma_X$, is non-singular.

The same assumptions hold for the
increments of $Y$ (in that case, we call $\Sigma_Y$ the covariance matrix of
$Y_1$).
\end{hyp}

For $\gb\in\R,N\in\N$ 
and for a fixed realization of
$Y$ 
we define a Gibbs transformation of the path measure $\Px$: this is
the polymer path measure $\Pnb$, absolutely continuous with
respect to $\Px$, given by
\begin{equation}
  \frac{\dd \Pnb}{\dd \Px} (X) = \frac{e^{\gb L_N(X,Y)}\; \ind_{\{\xyN\}}}{\Znb},
\end{equation}
where $L_N(X,Y)=\sum\limits_{n=1}^N \ind_{\{\xyn\}}$, 
and where
\begin{equation}
\Znb = \Ex[e^{\gb L_N(X,Y)}\; \ind_{\{\xyN\}}]
\end{equation}
is the partition function that normalizes $\Pnb$ to a probability.


The \textit{quenched} free energy  of the model is defined by
\begin{equation}
F(\gb):= \lim_{\Ninf} \frac{1}{N} \log \Znb = \lim_{\Ninf} \frac{1}{N}\Ey[ \log \Znb] 
\end{equation}
(the existence of the limit and the fact that it is $\bbP^Y$-almost surely 
constant is proven in \cite{BS08}).  We define also the $\ann$ partition
function $\Ey[\Znb]$, and the \textit{annealed} free energy:
\begin{equation}
\Fa(\gb):= \lim_{\Ninf} \frac{1}{N} \log \Ey[\Znb].
\end{equation}
We can compare the \textit{quenched} and \textit{annealed} free
energies, via the Jensen inequality:
\begin{equation}
 F(\gb)= \lim_{\Ninf} \frac{1}{N} \Ey[\log \Znb] \leq  \lim_{\Ninf} \frac{1}{N} \log \Ey[\Znb] = \Fa(\gb). \label{jensen free energy}
\end{equation}
The properties of $\Fa(\cdot)$ are well known (see the Remark~\ref{fctannealed}), and we have the existence of
critical points \cite{BS08}, for both \textit{quenched} and \textit{annealed} models, thanks to the convexity
and the monotonicity of the free energies with respect to  $\beta$:
\begin{definition}[Critical points]
 There exist $0\leq \gba \leq \gb_c$ depending on the laws of $X$ and $Y$ such that:
$\Fa(\gb)=0$ if $\gb\leq \gba$ and $\Fa(\gb)>0$ if $\gb>\gba$;
$F(\gb)=0$ if $\gb\leq \gb_c$ and $F(\gb)>0$ if $\gb>\gb_c$.
\end{definition}
The inequality $\gba \leq \gb_c$ comes from the inequality~(\ref{jensen free energy}).
\begin{rem} \rm 
As was remarked  in \cite{BS08}, 
the \textit{annealed} model is just the
 homogeneous pinning model \cite[Chapter 2]{GBbook} with partition function
$$\Ey[\Znb]=\Exy\left[\exp\left(\gb \sum_{n=1}^N \ind_{\{(X-Y)_n=0\}}\right)
\ind_{\{(X-Y)_N=0\}}
\right]$$
which describes the
 random walk $X-Y$ which receives the reward
$\beta$ each time it hits  $0$.
From the well-known results on the homogeneous pinning model  one sees
therefore that
\begin{itemize}
\item If $d=1$ or $d=2$, the $\ann$ critical point $\gba$ is
  zero because the random walk $X-Y$ is recurrent. 
	\item If $d\geq3$, the walk $X-Y$ is transient 
and as a consequence 
$$\gba = -\log\left[1-\Pxy\big( (X-Y)_n\ne 0\mbox{\;\;for every\;\;} n>0\big) \right]>0.$$
\end{itemize}
\label{fctannealed}
\end{rem}

\begin{rem}\rm
 As in the pinning model~\cite{GBbook}, the critical point $\gb_c$ marks the transition from a delocalized to a localized regime.
We observe that thanks to the convexity of the free energy,
\begin{equation}
 \partial_{\gb} F(\gb) = \lim_{N\to\infty} \Enb\left[ \frac{1}{N} \sum_{n=1}^N 
\ind_{\{\xyN\}} \right],
\end{equation}
almost surely in $Y$,
for every $\gb$ such that $F(\cdot)$ is differentiable at $\gb$. This
is the contact fraction between $X$ and $Y$.  When $\gb<\gb_c$, we
have $F(\gb)=0$, and the limit density of contact between $X$ and $Y$
is equal to~$0$: $\Enb\sum_{n=1}^N \ind_{\{\xyN\}}=o(N)$, and we are in the
delocalized regime. On the other hand, if $\gb>\gb_c$, we have
$F(\gb)>0$, and there is a positive density of contacts between $X$
and $Y$: we are in the localized regime.
\end{rem}

\subsection{Review of the known results}

The following is known about the question of the coincidence of
quenched and annealed critical points:
\begin{theorem} 
\cite{BS08}
Assume that $X$ and $Y$ are discrete time simple random walks on $\Z^d$.

If $d=1$ or $d=2$, the $quenched$ and $\ann$ critical points coincide: 
$\gbc=\gba=0.$

If $d\geq4$, the $quenched$ and $\ann$ critical points differ:
$\gbc>\gba>0.$
\label{quenched/annealed}
\end{theorem}
In dimension $d\ge 5$, the result was also proven (via a
very different method, and for more general random walks which include 
those of Assumption \ref{hypwright}) in an early
version of the paper \cite{BGdH08}.

\begin{rem}\rm
  The method and result of \cite{BS08} in dimensions $d=1,2$ can be
  easily extended beyond the simple random walk case (keeping zero mean 
and finite variance). On the other hand, in the case $d\ge4$ new ideas 
are needed to make the change-of-measure argument of \cite{BS08} work
for more general random walks.

Birkner and Sun gave also a similar result if $X$ and $Y$ are
continuous-time symmetric simple random walks on $\Z^d$, with jump rates $1$ and
$\rho\geq0$ respectively.  With definitions of (quenched and annealed)
free energy and critical points which are analogous to those of the 
discrete-time model, they proved:
\begin{theorem}\cite{BS08}
  In dimension $d=1$ and $d=2$, one has ${\gb}_c={\gb}_c^{ann}=0$.  In
  dimensions $d\geq 4$, one has $0<{\gb}_c^{ann}<{\gb}_c$ for each
  $\rho>0$.  Moreover, for $d=4$ and for each $\gd>0$, there exists
  $a_{\gd}>0$ such that ${\gb}_c-{\gb}_c^{ann}\geq a_{\gd}
  \rho^{1+\gd}$ for all $\rho\in[0,1]$.  For $d\geq 5$, there exists
  $a>0$ such that ${\gb}_c-{\gb}_c^{ann} \geq a \rho$ for
  all $\rho\in[0,1]$.
\end{theorem}
\end{rem}

Our main result completes this picture, resolving the open case of
the critical dimension $d=3$ (for simplicity,
 we deal only with the discrete-time model).
\begin{theorem}
\label{th:main}
Under the Assumption \ref{hypwright},
for $d=3$, we have
$\gb_c>\gba.$
\end{theorem}

We point out that the result holds also in the case where $X$ (or
$Y$) is a simple random walk, a case which a priori is excluded
by the aperiodicity condition of Assumption \ref{hypwright};
see the Remark~\ref{CLTsimple}.

Also, it is possible to modify our change-of-measure argument
to prove the non-coincidence of quenched and annealed critical points
in dimensions $d=4$ for the general walks of Assumption \ref{hypwright},
thereby extending 
the result of \cite{BS08}; see Section \ref{rem:d4} for a hint at the
necessary steps. 

{\bf Note} After this work was completed, M. Birkner and R. Sun
informed us that in \cite{BS09} they independently proved Theorem
\ref{th:main} for the continuous-time model.

\subsection{A renewal-type representation for $\Znb$}
From now on, we will assume that $d\ge 3$.

As discussed in \cite{BS08}, there is a way to represent the
partition function $Z_{N,Y}^\beta$ in terms of a renewal process
$\tau$; this rewriting makes the model look formally similar to
the random pinning model \cite{GBbook}. 

In order to introduce the representation of \cite{BS08}, we need a few
definitions. 
\begin{definition}
  We let
  \begin{enumerate}

  \item $\px_n(x)=\Px(X_n=x)$ and $\pxy_n(x)=\Pxy\big( (X-Y)_n=x \big)$;

  \item $\bP$ be the law of a recurrent renewal $\tau=\{\tau_0,\tau_1,\ldots\}$ with $\tau_0=0$, i.i.d. increments and inter-arrival
law given by 
\begin{eqnarray}
\label{K}
  K(n):=\bP(\tau_1=n)=\frac{\pxy_n(0)}{\Gxy} \mbox{\;\;where\;\;} \Gxy:=\sum_{n=1}^{\infty}\pxy_n(0)
\end{eqnarray}
(note that $\Gxy<\infty$ in dimension $d\ge 3$);

  \item $z'=(e^\gb-1)$ and $z=z'\,\Gxy$;

\item for $n\in\N$ and $x\in \Z^d$,
\begin{equation}
w(z,n,x)=z\frac{\px_{n}(x)}{\pxy_{n}(0)};
\label{defw}
\end{equation}

  \item $\check Z_{N,Y}^z:=\frac{z'}{1+z'}\Znb$.
 \end{enumerate}
\end{definition}

Then, via the binomial expansion of 
$e^{\gb L_N(X,Y)}=(1+z')^{L_N(X,Y)}$ one gets \cite{BS08}
\begin{eqnarray}
  \check Z_{N,Y}^z&=&
\sum_{m=1}^N \sum_{\tau_0=0<\tau_1<\ldots<\tau_m=N} \prod_{i=1}^m K(\tau_i-\tau_{i-1})w(z,\tau_i-\tau_{i-1},Y_{\tau_i}-Y_{\tau_{i-1}}) \label{decompZnzc}
\\\nonumber
&=&\E\left[W(z,\tau\cap\{0,\ldots,N\},Y)\ind_{N\in\tau}\right],
\end{eqnarray}
where we defined for any finite  increasing sequence $s=\{s_0,s_1,\ldots,s_l\}$
\begin{equation}
  W(z,s,Y) = \frac{\Ex\left[\left.\prod_{n=1}^{l} z \ind_{\{X_{s_n}=Y_{s_n}\}}\right|X_{s_0}=Y_{s_0}\right]}
  {\Exy\left[ \left.\prod_{n=1}^{l} \ind_{\{X_{s_n}=Y_{s_n\}}}\right|X_{s_0}=Y_{s_0}\right]}=\prod_{n=1}^{l} w(z,s_{n}-s_{n-1},Y_{s_n}-Y_{s_{n-1}}).
\label{Wtau}
\end{equation}

We remark that, taking the $\Ey-$expectation of the weights, we get
$${\Ey\left[w(z,\tau_i-\tau_{i-1},Y_{\tau_i}-Y_{\tau_{i-1}})\right]=z}.$$
Again, we see that the $\ann$ partition function is the partition function of a
homogeneous pinning model:
\begin{equation}
\check Z_{N,Y}^{z,ann}=\Ey[\check Z_{N,Y}^z]
=\bE\left[
z^{R_N}\ind_{\{N\in\tau\}}
\right],
\end{equation}
where we defined $R_N:=\left|\tau\cap\{1,\ldots,N\}\right|$.

Since the renewal $\tau$ is recurrent, the $\ann$ critical point is
$z_c^{ann}=1$.

\bigskip

In the following, we will often use the Local Limit Theorem for random walks. 
The following formulation can be extracted for instance from 
\cite[Theorem 3]{Stone} (recall that we assumed that the increments 
of both $X$ and $Y$ have finite exponential moments and non-singular 
covariance matrix):
\begin{proposition}[Local Limit Theorem]
\label{prop:localCLT}
Under the Assumption \ref{hypwright}, we get
\begin{equation}
  \bbP^X(X_n=x) = (1+o(1))\frac{1}{(2\pi n)^{d/2}(\det \gS_X)^{1/2}} 
\exp\left( -\frac{1}{2n} x\cdot\left(\gS_X^{-1} x\right)\right),
\end{equation}
where $o(1)\to 0$ as $n\to \infty$, uniformly for $\|x\|\leq n^{3/5}$.\\
Moreover, there exists a constant $c>0$ such that for all $x\in\Z^d$
\begin{equation}
 \bbP^X(X_n=x)\leq cn^{-d/2}.
\label{CLTbound}
\end{equation}
Similar statements hold for the walk $Y$.
\end{proposition}
(We use 
the notation $x\cdot y$ for the canonical
scalar product in $\R^d$.)

\smallskip

In particular, from Proposition \ref{prop:localCLT} and the definition
of $K(\cdot)$ in \eqref{K}, we get $K(n)\sim c_K n^{-d/2}$ as
$n\to\infty$, for some positive $c_K$. As a consequence, we 
get from  \cite[Th. B]{Doney}
that
\begin{eqnarray}
  \label{eq:doney}
  \bP(n\in\tau)\stackrel{n\to\infty}\sim \frac1{2\pi c_K\sqrt n}.
\end{eqnarray}

\begin{rem}\rm
\label{CLTsimple}
In Proposition~\ref{prop:localCLT}, we supposed that the walk $X$ is
aperiodic, which is not the case for the simple random walk. 
If $X$ is the symmetric simple random walk on 
 $\Z^d$, then \cite[Prop. 1.2.5]{Lawl}
\begin{equation}
  \bbP^X(X_n=x) = (1+o(1))\ind_{\{n\leftrightarrow x\}}\frac{2}{(2\pi n)^{d/2}(\det \gS_X)^{1/2}} 
\exp\left( -\frac{1}{2n} x\cdot\left(\gS_X^{-1} x\right)\right),
\end{equation}
where $o(1)\to 0$ as $n\to \infty$, uniformly for $\|x\|\leq n^{3/5}$,
and where $n\leftrightarrow x$ means that $n$ and $x$ have the same
parity (so that $x$ is a possible value for $X_n$). Of course, in this
case $\gS_X$ is just $1/d$ times the identity matrix.
The
statement~\eqref{CLTbound} also holds.

Via this remark, one can adapt all the computations of the following
sections, which are based on Proposition~\ref{prop:localCLT}, to the
case where $X$ (or $Y$) is a simple random walk. For 
simplicity of exposition, we give the proof of Theorem \ref{th:main}
only in the aperiodic case.
\end{rem}

\section{Main result: the dimension $d=3$}
\label{section d=3}

With the definition $\check F(z):=\lim_{\Ninf}\frac{1}{N}\log \check 
Z_{N,Y}^z$, to prove Theorem \ref{th:main} it is sufficient to
show that $\check F(z)=0$  for some $z>1$.

\subsection{The coarse-graining procedure and the fractional moment method}

We consider without loss of generality a system of size proportional to $L=\frac{1}{z-1}$ (the coarse-graining length), that is $N=mL$, with $m\in \N$.
Then, for $\cI \subset \{1,\ldots,m\}$, we define
\begin{equation}
  \ZI := \E\big[ W(z,\tau\cap\{0,\ldots,N\},Y)\ind_{N\in\tau} \ind_{ E_{\cI}}(\tau)\big] ,
\end{equation}
where $E_{\cI}$ is the event that the renewal $\tau$ intersects the blocks $(B_i)_{i\in\cI}$ and only these blocks over $\{1,\ldots,N\}$, $B_i$ being the $i^{th}$ block of size $L$:
\begin{equation}
B_i:=\{(i-1)L+1,\ldots,i L\}.
\end{equation}
Since the events $E_\cI$ are disjoint, we can write
\begin{equation}
\check Z^z_{N,Y} := \sum_{\cI\subset \{1,\ldots,m\}} \ZI.
\label{decoupZi}
\end{equation}
Note that $\ZI=0$ if $m\notin \cI$. We can therefore assume $m\in\cI$.
If we denote $\cI=\{i_1,i_2,\ldots,i_l\}$ ($l=\left|\cI\right|$),
$i_1<\ldots<i_l$, $i_l=m$, we can express $\ZI$ in the following way:

\begin{eqnarray}
\ZI & := & \sumtwo{a_1,b_1 \in B_{i_1}}{a_1\leq b_1} \sumtwo{a_2,b_2 \in B_{i_2}}{a_2\leq b_2}  \ldots \sum_{a_l \in B_{i_l}} K(a_1)w(z,a_1,Y_{a_1}) Z_{a_1,b_1}^z  \label{coarsedecomp}\\
& & \indent  \indent \indent \indent \indent \ldots K(a_l-b_{l-1})w(z,a_l-b_{l-1},Y_{a_l}-Y_{b_{l-1}}) Z_{a_l,N}^z,\nonumber 
\end{eqnarray}
where
\begin{equation}
Z_{j,k}^z := \E\big[  W(z,\tau\cap\{j,\ldots,k\},Y)\ind_{k\in\tau}\left|j\in\tau\right. \big]
\label{Zjk}
\end{equation}
is the  partition function between $j$ and $k$.

\begin{figure}[htbp]
\centerline{
\psfrag{0}{$0$}
\psfrag{L}{$L$}
\psfrag{2L}{$2L$}
\psfrag{3L}{$3L$}
\psfrag{4L}{$4L$}
\psfrag{5L}{$5L$}
\psfrag{6L}{$6L$}
\psfrag{7L}{$7L$}
\psfrag{8L}{$8L=N$}
\psfrag{d1}{$a_1$}
\psfrag{d2}{$a_2$}
\psfrag{d3}{$a_3$}
\psfrag{d4}{$a_4$}
\psfrag{f1}{$b_1$}
\psfrag{f2}{$b_2$}
\psfrag{f3}{$b_3$}
\psfrag{f4}{$b_4=N$}
\psfig{file=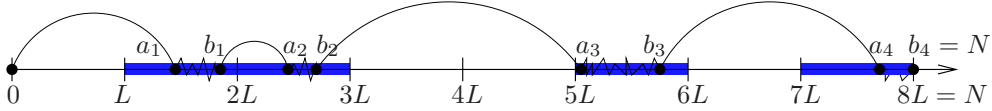,width=5in} }
   \caption{ The coarse-graining procedure. Here $N=8L$ (the system is cut into $8$ blocks),
and $\cI = \{ 2,3,6,8 \}$ (the gray zones) are the blocks where the contacts occur, and where the change of measure procedure
of the Section~\ref{chgtmes} acts.}
\label{figcoarse}
\end{figure}

Moreover, thanks to the Local Limit Theorem (Proposition~\ref{prop:localCLT}), one can note that there
exists a constant $c>0$ independent of the realization of $Y$ such that, if one takes
$z\leq 2$ (we will take $z$ close to $1$ anyway), one has
$$w(z,\tau_{i}-\tau_{i-1},Y_{\tau_i}-Y_{\tau_{i-1}})=z\frac{\px_{\tau_i-\tau_{i-1}}(Y_{\tau_i}-Y_{\tau_{i-1}})}{\pxy_{\tau_i-\tau_{i-1}}(0)} \le c.
$$
So, the decomposition~(\ref{coarsedecomp}) gives
\begin{equation}
\ZI \leq c^{\left|\cI\right|} \sumtwo{a_1,b_1 \in B_{i_1}}{a_1\leq b_1} \sumtwo{a_2,b_2 \in B_{i_2}}{a_2\leq b_2}  \ldots \sum_{a_l \in B_{i_l}}
 K(a_1) Z_{a_1,b_1}^z K(a_2-b_1) Z_{a_2,b_2}^z \ldots K(a_l-b_{l-1})Z_{a_l,N}^z.
\label{ineqcoarse}
\end{equation}

We now eliminate the dependence on $z$ in the
inequality~(\ref{ineqcoarse}). This is possible thanks to the choice
$L=\frac{1}{z-1}$. As each $Z_{a_i,b_i}^z$ is the partition
function of a system of size smaller than $L$, we get
$W(z,\tau\cap\{a_i,\ldots,b_i\},Y)\leq z^L
W(z=1,\tau\cap\{a_i,\ldots,b_i\},Y)$ (recall the definition \ref{Wtau}).
 But with the choice $L=\frac{1}{z-1}$, the factor $z^L$ is
bounded by a constant $c$, and thanks to the equation~(\ref{Zjk}),
we finally get
\begin{equation}
Z_{a_i,b_i}^z\leq c Z_{a_i,b_i}^{z=1}.
\label{dependz}
\end{equation}

{\bf Notational warning}: in the following, $c,c'$, etc. will denote positive constants,
whose value may change from line to line. 

\smallskip

We note $Z_{a_i,b_i}:=Z_{a_i,b_i}^{z=1}$ and
$W(\tau,Y):=W(z=1,\tau,Y)$. Plugging this in the
inequality~(\ref{ineqcoarse}), we finally get
\begin{equation}
\ZI \leq c'^{\left|\cI\right|} \sumtwo{a_1,b_1 \in B_{i_1}}{a_1\leq b_1} \sumtwo{a_2,b_2 \in B_{i_2}}{a_2\leq b_2}  \ldots \sum_{a_l \in B_{i_l}}
 K(a_1) Z_{a_1,b_1} K(a_2-b_1) Z_{a_2,b_2} \ldots K(a_l-b_{l-1})Z_{a_l,N},
\label{coarseineq}
\end{equation}
where there is no dependence on $z$ anymore.

The fractional moment method starts from the observation  that for any
$\gamma\ne0$
\begin{equation}
  \check F(z) = 
\lim_{N\to\infty} \frac{1}{\gga N} 
\Ey\left[\log \left(\check Z^z_{N,Y}\right)^{\gga}\right] \leq \liminf_{\Ninf} 
\frac1{N\gamma}\log \Ey\left[\left(\check Z^z_{N,Y}\right)^{\gga}\right].
\label{momentfractionnaire}
\end{equation}
Let us fix a value of $\gga\in(0,1)$ (as in \cite{GLT09}, we will
choose $\gga=6/7$, but we will keep writing it as $\gga$ to simplify
the reading). Using the inequality $\left(\sum a_n\right)^{\gga}\leq
\sum a_n^{\gga}$ (which is valid for $a_i\geq 0$), and combining with
the decomposition~\eqref{decoupZi}, we get
\begin{equation}
\Ey\left[\left(\check Z^z_{N,Y}\right)^{\gga}\right] \leq  \sum_{\cI\subset \{1,\ldots,m\}} \Ey\left[\left(\ZI\right)^{\gga}\right].
\label{momfrac}
\end{equation}

Thanks to~(\ref{momentfractionnaire}) we only have to prove that, for some $z>1$, $\limsup_{\Ninf} \Ey\left[\left(
\check Z_{N,Y}^z\right)^{\gga}\right]<\infty$.\\
We deal with the term $\Ey\left[(\ZI)^{\gga}\right]$ via a change of measure procedure.

\subsection{The change of measure procedure}
\label{chgtmes}
The idea is to change the measure $\Py$ on each block whose index
belongs to $\cI$, keeping each block independent of the others.  We
replace, for fixed $\cI$, the measure $\Py(\dd Y)$ with
$g_{\cI}(Y)\Py(\dd Y)$, where the function $g_{\cI}(Y)$ will have the effect of creating
 long range positive correlations between the increments
of $Y$, inside each block separately.
Then, thanks to the H\"older inequality, we can write
\begin{eqnarray}
\Ey\left[\left(\ZI\right)^{\gga}\right]=\Ey\left[
\frac{g_{\cI}(Y)^\gamma}{g_{\cI}(Y)^\gamma}\left(\ZI\right)^{\gga}\right] & \leq & \Ey\left[g_{\cI}(Y)^{-\frac{\gga}{1-\gga}}\right]^{1-\gga}\Ey\left[g_{\cI}(Y)\ZI\right]^{\gga}.
\label{holder}
\end{eqnarray}

In the following, we will denote $\gD_i=Y_{i}-Y_{i-1}$ the $i^{th}$
increment of $Y$. 
Let us introduce, for $K>0$ to
be chosen, the following ``change of measure'':
\begin{eqnarray}
  \label{eq:change1block}
  g_{\cI}(Y)= 
\prod_{k\in\cI} e^{-F_k(Y)\ind_{F_k(Y)\geq0}}
\equiv \prod_{k\in\cI}g_k(Y),
\end{eqnarray}
where
\begin{equation}
F_k(Y) = -\sum\limits_{i,j\in B_k}{M_{ij} \gD_i\cdot\gD_j},
\label{fun chgtmes}
\end{equation}
and
\begin{equation}
\begin{cases}
M_{ij}=\frac{c_M}{\sqrt{L\log L}}\frac{1}{\sqrt{\left|j-i\right|}} & \text{if }\  i\neq j \\
M_{ii}=0. & 
\end{cases}
\label{defM}
\end{equation}
The constant  $c_M$ will be chosen in a moment. We note that $F_k$ only
depends on the increments of $Y$ in the block labeled $k$. 
\paragraph{Let us deal with the first factor of~(\ref{holder}):}
\begin{equation}
\Ey\left[g_{\cI}(Y)^{-\frac{\gga}{1-\gga}}\right] = \prod_{k\in\cI}
 \Ey\left[g_k(Y)^{-\frac{\gga}{1-\gga}}\right] = 
\Ey\left[\exp\left(\frac\gga{1-\gga}F_1(Y)\ind_{F_1(Y)\ge0}\right)
\right]^{\left|\cI\right|}.
\end{equation}
We use the following lemma to choose $c_M$:
\begin{lemma}
There exists a constant $c>0$, such that if $\norm{M}^2:=\sum_{i,j\in B_1}M_{ij}^2<c$, then
\begin{eqnarray}
  \label{eq:borneM}
\Ey\left[\exp \left(\frac{\gga}{1-\gga} 
F_1(Y)\ind_{F_1(Y)\ge0}
\right)\right] \leq 4.
\end{eqnarray}
\label{borneM}
\end{lemma}
It is immediate to check that if $c_M$ is taken small enough, then $\norm{M}^2<c$, and
the first factor in \eqref{holder} is bounded by $4^{(1-\gga)\left|\cI\right|}\leq
3^{\gga \left|\cI\right|}$. The inequality~(\ref{holder}) finally
gives
\begin{equation}
\Ey\left[\left(\ZI\right)^{\gga}\right] \leq 3^{\gga|\cI|} \Ey\left[g_{\cI}(Y)\ZI\right]^{\gga}.
\label{ZIgamma}
\end{equation}
{\sl Proof of Lemma \ref{borneM}}.
We use the following theorem, proved in \cite{Wri}, which gives a bound for the tail probability for quadratic forms in independent random variables.
\begin{theorem}
\label{th:wright}
  Let $(X_i)_{i\geq 1}$ be a sequence of independent real-valued random variables,
  with zero means and which verify $P(|X_i| \geq x ) \leq M
  \int_x^{\infty} exp(-h t^2) \dd t $ for all $x\geq0$, where $M$
  and $h$ are positive constants. Let $A$ be a real symmetric
  matrix such that $\norm{A}^2=\sum_{ij} a_{ij}^2 <\infty$,
  let $\rho(A)$ be the norm of $A$ considered as an operator in $l_2$ and set $S=\sum a_{ij}\left(X_i X_j- E(X_i X_j)\right)$.\\
  Then, there exist two positive constants $C_1$ and $C_2$ (which depend only on
  $M$ and $h$) such that, for every $u\geq 0$
$$ P\left(S\geq u\right) \leq  \exp \left\{ -\min \left(\frac{C_1 u}{\rho(A)}, \frac{C_2 u^2}{\norm{A}^2}\right) \right\}.$$
\end{theorem}
It is well known that the Hilbert-Schmidt norm of a matrix dominates its spectral radius, $\rho(A)\leq \norm{A}$. 
Then, it follows that 
for every $u\geq (C_1/C_2)\|A\|$ one has
\begin{equation}
P\left(S\geq u\right) \leq  \exp \left(-\frac{C_1 u}{\norm{A}}\right).
\end{equation}

In order to prove \eqref{eq:borneM}, we introduce $\gD_i =(\gD_i^{(1)},\gD_i^{(2)},\gD_i^{(3)})$, the components of
$\gD_i$ 
and
$$S=-\sum_{e=1}^{3}\sum_{i,j\in B_1} \frac{\gga}{1-\gga} {M_{ij} \gD^{(e)}_i \gD^{(e)}_j} = 
\sum_{e=1}^{3} S^{(e)}.$$ We can then apply  Theorem \ref{th:wright} for $S^{(e)}$
with $X_i = \gD_i^{(e)}$ (recall the assumption \ref{hypwright}) and
$a_{ij}=-\frac{\gga}{1-\gga}M_{ij}$ for $i,j\leq L$ (we take
$a_{ij}=0$ if $i>L$ or $j>L$, and recall $M_{ii}=0$).  
Let us choose $c_M$ sufficiently small so that $(C_1/C_2)\|A\|\le 1/3$.
For
any $u\ge 1$, we get
\begin{equation}
\Py( S\geq u )  \leq \sum_{k=1}^{3} \Py(S^{(k)}\geq u/3)  \leq  3 e^{-\frac{C_1(1-\gamma)}{3\gamma\norm{M}}u}.
\end{equation}
Then,
\begin{eqnarray*}
\Ey\left[e^{S\ind_{S\ge0}}\right]  \le 1+\Ey\left[e^S\right]\le
1+e+
\int_{e}^{+\infty} \Py(S \geq \log u) \dd u 
  \leq 1+ e+ 3 \int_{e}^{+\infty} u^{-\frac{C_1(1-\gamma)}{3\gamma\norm{M}}}\dd u.
\end{eqnarray*}
Choosing $c_M$ sufficiently small, the right-hand side can be made
as close to $1+e<4$ as wished.
\qed

\medskip

We are left with the estimation of $\Ey\left[g_{\cI}(Y)\ZI\right]$.
We set $\PI := \P\left(\EI, N\in\tau \right)$, that is the probability for $\tau$ to visit the blocks $(B_i)_{i\in\cI}$ and only these ones, and to 
visit also $N$. We now use the following two statements.

\begin{proposition}
\label{th:core}
For any $\eta>0$, there exists $z>1$ sufficiently close to $1$ (or $L$ sufficiently big, since $L=(z-1)^{-1}$) such that for every $\cI\subset \{1,\ldots,m\}$ with $m\in\cI$, we have
\begin{equation}
\Ey\left[g_{\cI}(Y)\ZI\right] \leq \eta^{\left|\cI\right|} \PI.
\end{equation}
\label{EgZI}
\end{proposition}
Proposition \ref{th:core} is the core of the paper and is proven in
the next section.

\begin{figure}[htbp]
\centerline{
\psfrag{0}{$0$}
\psfrag{1}{$1$}
\psfrag{2}{$2$}
\psfrag{3}{$3$}
\psfrag{4}{$4$}
\psfrag{5}{$5$}
\psfrag{6}{$6$}
\psfrag{7}{$7$}
\psfrag{8}{$8=m$}
\psfig{file=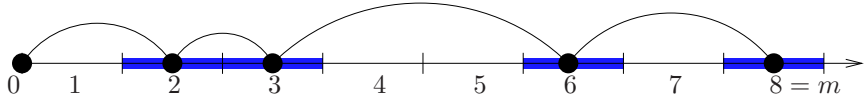,width=4.5in} }
   \caption{  The underlying renewal $\tilde{\tau}$ is a
subset of the set of blocks $(B_i)_{1\leq i\leq n}$ (i.e the blocks are
reinterpreted as points) and the inter-arrival distribution is $\tilde{K}(n)= 1 /\left(\tilde c n^{6/5}\right)$.}
\label{coarsepinning}
\end{figure}

\begin{lemma}\cite[Lemma 2.4]{GLT09}
There exist three constants $C_1=C_1(L)$, $C_2$ and $L_0$  such that (with $i_0:=0$)
\begin{equation}
\PI \leq C_1 C_2^{\left|\cI\right|} \prod_{j=1}^{\left|\cI\right|} \frac{1}{(i_j-i_{j-1})^{7/5}}
\end{equation}
for $L\geq L_0$ and for every $\cI\in\{1,\ldots,m\}$.
\label{PI}
\end{lemma}
Thanks to these two statements and combining with the inequalities~(\ref{momfrac}) and~(\ref{ZIgamma}), we get
\begin{equation}
  \Ey\left[\left(\Znzc\right)^{\gga}\right]\leq \sum_{\cI\subset \{1,\ldots,m\}} 
\Ey\left[\left(\ZI\right)^{\gga}\right] \leq C_1^{\gga} \sum_{\cI\subset \{1,\ldots,m\}}  
\prod_{j=1}^{\left|\cI\right|}\frac{(3C_2\eta)^{\gga} }{(i_j-i_{j-1})^{7\gga/5}}.
\label{EZgga}
\end{equation}
Since $7\gga/5 = 6/5 >1$, we can set
\begin{equation}
\tilde{K}(n) = \frac{1}{\tilde c n^{6/5}}, \indent \text{where } \ \tilde{c} = \sum_{i=1}^{+\infty} i^{-6/5} <+ \infty,
\end{equation}
and $\tilde{K}(\cdot)$ is the inter-arrival probability of some recurrent renewal
$\tilde\tau$.  We can therefore interpret
the right-hand side of~(\ref{EZgga}) as a partition function of a homogeneous pinning
model of size $m$ (see Figure~\ref{coarsepinning}), with the underlying renewal $\tilde{\tau}$, and
with pinning parameter $\log[\tilde c(3C_2\eta)^{\gga}]$:
\begin{equation}
\Ey\left[\left(\check Z^z_{N,Y}\right)^{\gga}\right] \leq C_1^{\gga} \E_{\tilde{\tau}} \left[\left(\tilde c(3C_2\eta)^{\gga}\right)^{\left|\tilde{\tau}\cap \{1,\ldots,m\}\right|}\right].
\label{EZgga2}
\end{equation}

Thanks to  Proposition~\ref{EgZI}, we can take $\eta$ arbitrary small. Let us fix
$\eta := 1/((4C_2)\tilde{c}^{1/\gga})$.
Then,
\begin{equation}
\Ey\left[\left(\check Z^z_{N,Y}\right)^{\gga}\right] \leq C_1^{\gga}
\end{equation}
for every $N$. This implies, thanks to~(\ref{momentfractionnaire}), that $\check
F(z)=0$, and we are done. \qed

\begin{rem}\rm
  The coarse-graining procedure reduced the proof of delocalization to
  the proof of Proposition~\ref{EgZI}. Thanks to the
  inequality~(\ref{coarseineq}), one has to estimate the expectation,
  with respect to the $g_{\cI}(Y)-$modified measure, of the partition
  functions $Z_{a_i,b_i}$ in each visited block. We will show (this is
  Lemma~\ref{lemEgZ}) that the expectation with respect to this
  modified measure of $Z_{a_i,b_i}/\P(b_i-a_i\in\tau)$ can be
  arbitrarily small if $L$ is large, and if $b_i-a_i$ is of the order of
  $L$. If $b_i-a_i$ is much smaller, we can deal with this term via elementary
bounds.
\label{sketchprop}
\end{rem}

\section{Proof of the Proposition~\ref{EgZI}}
\label{section prop EgZI}

As pointed out in Remark~\ref{sketchprop}, Proposition~\ref{EgZI}
relies on the following key lemma: 
\begin{lemma}
For every $\gep$ and $\gd>0$, there exists $L>0$ such that
\begin{equation}
\Ey\left[g_1(Y) Z_{a,b}\right] \leq \gd \P(b-a\in\tau)
\end{equation}
for every $a\leq b$ in $B_1$ such that $b-a\geq \gep L$.
\label{lemEgZ}
\end{lemma}

Given this lemma, the proof of Proposition~\ref{EgZI} is very similar
to the proof of \cite[Proposition 2.3]{GLT09}, so we will sketch only
a few steps.
The
 inequality~(\ref{coarseineq}) gives us
\begin{eqnarray}
\lefteqn{\Ey\left[g_{\cI}(Y)\ZI\right]  } \nonumber\\
 & \leq  & c^{\left|\cI\right|} \sumtwo{a_1,b_1 \in B_{i_1}}{a_1\leq b_1} \sumtwo{a_2,b_2 \in B_{i_2}}{a_2\leq b_2}  \ldots \sum_{a_l \in B_{i_l}}
 K(a_1) \Ey\left[g_{i_1}(Y)Z_{a_1,b_1}\right] K(a_2-b_1) \Ey\left[g_{i_2}(Y)Z_{a_2,b_2}\right] \ldots \nonumber\\
  & & \indent\indent\indent \ldots K(a_l-b_{l-1})\Ey\left[g_{i_l}(Y)
Z_{a_l,N}\right] \nonumber
 \nonumber\\
 & = & c^{\left|\cI\right|} \sumtwo{a_1,b_1 \in B_{i_1}}{a_1\leq b_1} \sumtwo{a_2,b_2 \in B_{i_2}}{a_2\leq b_2}  \ldots \sum_{a_l \in B_{i_l}}
 K(a_1) \Ey\left[g_1(Y)Z_{a_1-L(i_1-1),b_1-L(i_1-1)}\right] K(a_2-b_1) \ldots \label{interm ineqEgZI}\\
  & & \indent\indent\indent \ldots K(a_l-b_{l-1})
\Ey\left[g_{1}(Y)Z_{a_l-L(m-1),N-L(m-1)}\right]. \nonumber
\end{eqnarray}
The terms with $b_i-a_i\ge \gep L$ are dealt with via  Lemma~\ref{lemEgZ}, while
for the remaining ones we just observe that
$\bbE^Y[g_1(Y)Z_{a,b}]\le \bP(b-a\in\tau)$ since
$g_1(Y)\le 1$.
One has then
\begin{eqnarray}
  \Ey\left[g_{\cI}(Y)\ZI\right] & \leq & c^{\left|\cI\right|} \sumtwo{a_1,b_1 \in B_{i_1}}{a_1\leq b_1} \sumtwo{a_2,b_2 \in B_{i_2}}{a_2\leq b_2}  \ldots \sum_{a_l \in B_{i_l}} K(a_1) \left(\gd+\ind_{\{b_1-a_1 \leq \gep L\}}\right)\P(b_1-a_1\in\tau)
  \nonumber \\
  & & \indent \indent \indent \ldots K(a_l-b_{l-1}) \left(\gd+\ind_{\{N-a_l \leq \gep L\}}\right)\P(N-a_l\in\tau). \label{ineqEgZI}
\end{eqnarray}
From this point on, the proof of Theorem \ref{th:core} is identical to the proof of
Proposition 2.3 in \cite{GLT09} (one needs of course to choose 
$\gep=\gep(\eta)$ and $\delta=\delta(\eta)$ sufficiently small). \qed

\subsection{Proof of Lemma~\ref{lemEgZ}}

Let us fix $a,b$ in $B_1$, such that $b-a\ge \gep L$. The small constants $\delta$ and $\gep$ are also fixed.
We recall that for a fixed configuration of $\tau$ such that
$a,b\in\tau$, we have $\Ey\big[ W(\tau\cap\{a,\ldots,b\},Y)\big] =1$ because $z=1$. We can therefore introduce the
probability measure (always for fixed $\tau$) 
\begin{equation}
\dd \Pct (Y) = W(\tau\cap\{a,\ldots,b\},Y) \dd \Py(Y)
\label{defPct}
\end{equation}
where we do not indicate the dependence on $a$ and $b$.
Let us note for later convenience that, in the particular case  $a=0$, the definition \eqref{Wtau} of $W$ implies that
for any function $f(Y)$
\begin{eqnarray}
\label{obsPct}
  \bbE_\tau[f(Y)]=\bbE^X\bbE^Y\big[ f(Y)|X_i=Y_i\;\forall i\in \tau\cap\{1,\ldots,b\}\big] .
\end{eqnarray}

\smallskip

Then, with the definition~(\ref{Zjk}) of $Z_{a,b}:=Z_{a,b}^{z=1}$, we get
\begin{equation}
\Ey\left[g_1(Y)Z_{a,b}\right] = 
\Ey \E \big[ g_1(Y) W(\tau\cap \{a,\ldots,b\},Y)\ind_{b\in\tau}\left|a\in \tau\right.\big] = \Edf \Ectdf[g_1(Y)] \P(b-a\in\tau),
\end{equation}
where $\hat \bP(\cdot):=\bP(\cdot|a,b\in\tau)$, 
and therefore we have to show that $\Edf\Ectdf[g_1(Y)]\leq \gd$.

With the definition~\eqref{eq:change1block} of $g_1(Y)$, we get that for any $K$
\begin{eqnarray}
  \Edf\Ect[g_1(Y)]    
 \leq  e^{-K}+\Edf \bbP_\tau\left(F_1 < K \right). \label{step1}
\end{eqnarray}
If we choose $K$ big enough, the first term is
smaller than $\gd/3$. We now use two lemmas to deal with the second
term. The idea is to first prove that $\Ect[F_1]$ is big with a 
$\Pdf-$probability close to $1$, and then that its variance is not too large.
\begin{lemma} 
  Let $a,b\in B_1$ satisfy $b-a\ge \gep L$. 
Then, for every
   $\gz>0$, one can find  two constants $g>0$ and $L_0>0$, such that
\begin{equation}
\label{eq:Fbig}
\Pdf\left(\Ect[F_1] \le g\sqrt{\log L}\right) \le \gz,
\end{equation}
for every $L\geq L_0$.
\label{EctF big}
\end{lemma}
Choose $\gz =\gd/3$ and fix $g>0$ such that \eqref{eq:Fbig} holds
for every $L$ sufficiently large. If $2K = g\sqrt{\log L} $ (and
therefore we can make $e^{-K}$ small enough by choosing $L$ large), we
get that
\begin{eqnarray}
  \Edf \bbP_\tau\left(F_1 < K \right) & \leq & \EL\bbP_\tau\big[F_1-\Ect[F_1] \leq-K\big] + \PL\left(\Ect[F_1] 
\leq 2K\right) \\
  & \leq & \frac{1}{K^2} \EL\Ect\left[\left(F_1-\Ect[F_1]\right)^2\right] + \gd/3 .
\end{eqnarray}
Putting this together with~(\ref{step1}) and with our choice of $K$, we have
\begin{equation}
\EL\Ect[g_1(Y)] \leq 2\gd/3 + \frac{4}{g^2 \log L} \EL\Ect\left[\left(F_1-\Ect[F_1]\right)^2\right]
\end{equation}
for $L\geq L_0$. Then we just have to prove that $\EL\Ect\left[\left(F_1-\Ect[F_1]\right)^2\right]=o(\log L)$. 
Indeed,
\begin{lemma} Let $a,b\in B_1$ satisfy $b-a\geq \gep L$.
Then there exists some constant $c>0$ such that
\begin{equation}
\EL\Ect\left[\left(F_1-\Ect[F_1]\right)^2\right] \leq c \left(\log L\right)^{3/4}
\end{equation}
for every $L>1$.
\label{variance}
\end{lemma}
We finally get that
\begin{equation}
\EL\Ect[g_1(Y)] \leq 2\gd/3 + c (\log L)^{-1/4},
\end{equation}
and there exists a constant $L_1>0$ such that for $L>L_1$
\begin{equation}
\EL\Ect[g_1(Y)] \leq \gd.
\end{equation}
\qed

\subsection{Proof of Lemma~\ref{EctF big}}
\label{lemma EctF big}
Up to now, the proof of Theorem \ref{th:main} is quite similar to the proof of the main result in \cite{GLT09}. 
Starting from the present section, instead, new ideas and technical results are needed.

Let us fix a realization of $\tau$ such that $a,b\in\tau$ (so that it has
a non-zero probability under $\hat \bP$) and let us 
note $\tau \cap\{a,\ldots b\}=\{\tau_{R_a}=a,\tau_{R_a+1},\ldots,\tau_{R_b}=b\}$
(recall that $R_n=|\tau\cap\{1,\ldots,n\}|$).
We observe (just go back to
the definition of $\bbP_\tau$) that, if $f$ is a function of the
increments of $Y$ in $\{\tau_{n-1}+1,\ldots,\tau_n\}$, $g$ of the
increments in $\{\tau_{m-1}+1,\ldots,\tau_m\}$ with $R_a<n\ne m\le R_b$, and if $h$ is
a function of the increments of $Y$ not in $\{a+1,\ldots,b\}$ then
\begin{eqnarray}
\label{dif intervals}
  \lefteqn{\bbE_\tau\big[ f\big(\{\Delta_i\}_{i\in \{\tau_{n-1}+1,\ldots,\tau_n\}}\big)
    g\big(\{\Delta_i\}_{i\in \{\tau_{m-1}+1,\ldots,\tau_m\}}\big)
    h\big(\{\Delta_i\}_{i\notin \{a+1,\ldots,b\}}\big)\big] }\\\nonumber
  &=&
  \bbE_\tau\big[ f\big(\{\Delta_i\}_{i\in \{\tau_{n-1}+1,\ldots,\tau_n\}}\big)\big]
  \bbE_\tau \big[ g\big(\{\Delta_i\}_{i\in \{\tau_{m-1}+1,\ldots,\tau_m\}}\big)\big]
  \bbE^Y \big[ h\big(\{\Delta_i\}_{i\notin \{a+1,\ldots,b\}}\big)\big],
\end{eqnarray}
and that
\begin{eqnarray}
\nonumber
  \bbE_\tau\big[ f\big(\{\Delta_i\}_{i\in \{\tau_{n-1}+1,\ldots,\tau_n\}}\big)\big]
 = \bbE^X\bbE^Y\big[f\big(\{\Delta_{i}\}_{i\in \{\tau_{n-1}+1,\ldots,\tau_n\}}\big)|
    X_{\tau_{n-1}}  =  Y_{\tau_{n-1}},X_{\tau_{n}} =   Y_{\tau_{n}}\big] \\
 =\bbE^X\bbE^Y\big[f\big(\{\Delta_{i-\tau_{n-1}}\}_{i\in \{\tau_{n-1}+1,\ldots,\tau_n\}}\big)|
    X_{\tau_n-\tau_{n-1}}  =  Y_{\tau_n-\tau_{n-1}}\big] .
\label{ident}
\end{eqnarray}

We want to estimate $\Ect[F_1]$: since the increments $\Delta_i$ for
$i\in B_1\setminus \{a+1,\ldots,b\}$ are i.i.d. and centered (like under $\bbP^Y$), we have
\begin{equation}
\label{F_1}
\Ect[F_1] := \sum_{i,j=a+1}^{b} M_{ij} \Ect[-\gD_i\cdot\gD_j].
\end{equation}
Via a time translation, one can always assume that $a=0$
and we do so from now on.

The key point is the following
 \begin{lemma}
\label{th:Ar}
   \begin{enumerate}
   \item If there exists $1\le n\le R_b$ such that
     $i,j\in\{\tau_{n-1}+1,\ldots,\tau_n\}$, then
\begin{eqnarray}
  \Ect[-\gD_i\cdot\gD_j]  = A(r)\stackrel {r\to\infty}\sim
\frac{C_{X,Y}}r
\label{A(r)}
\end{eqnarray}
 where $r=\tau_n-\tau_{n-1}$ (in particular, note that the 
expectation depends only on $r$) and $C_{X,Y}$ is a positive constant
which depends on $\bbP^X,\bbP^Y$;

\item otherwise, $\Ect[-\gD_i\cdot\gD_j] =  0$.
   \end{enumerate}
 \end{lemma}

{\sl Proof of Lemma \ref{th:Ar}}
{\bf Case (2).} Assume that $\tau_{n-1}<i\le \tau_n$ and $\tau_{m-1}<j\le \tau_m$
with $n\ne m$. 
Thanks to \eqref{dif intervals}-\eqref{ident} we have
that 
\begin{equation}
\Ect[\gD_i \cdot\gD_j]= \bbE^X\bbE^Y[\Delta_i|X_{\tau_{n-1}}=Y_{\tau_{n-1}},
X_{\tau_{n}}=Y_{\tau_{n}}]\cdot
\bbE^X\bbE^Y[\Delta_j|X_{\tau_{m-1}}=Y_{\tau_{m-1}},
X_{\tau_{m}}=Y_{\tau_{m}}]
\end{equation}
and both factors are immediately seen to be zero, since 
the laws of $X$ and $Y$ are assumed to be symmetric.

{\bf Case (1).} Without loss of generality, assume that $n=1$, 
so we only have to compute 
\begin{equation}
\Ey\Ex\left[\gD_i\cdot \gD_j\left|X_r=Y_r\right.\right].
\end{equation}
where $r=\tau_1$.  Let us fix $x\in\Z^3$, and 
denote $\Eyrx[\cdot]=\Ey[\cdot\left|Y_r=x\right.]$.
\begin{eqnarray*}
  \Ey[\gD_i\cdot\gD_j\left|Y_r=x\right.] & = & \Eyrx\left[\gD_i
\cdot\Eyrx\left[\gD_j\left|\gD_i\right.\right]\right] \\
  & = & \Eyrx\left[\gD_i\cdot\frac{x-\gD_i}{r-1}\right] = \frac{x}{r-1}
\cdot\Eyrx\left[\gD_i\right] - \frac{1}{r-1}\Eyrx\left[\norm{\gD_i}^2\right] \\
  & = & \frac{1}{r-1}\left(\frac{\norm{x}^2}{r}-\Eyrx\left[\norm{\gD_1}^2\right] \right),
\end{eqnarray*}
where we used the fact that under $\Pyrx$ the law of the increments
$\{\gD_i\}_{i\le r}$ is exchangeable.
Then, we get
\begin{eqnarray}
 \lefteqn{ \Ect[\gD_i\cdot\gD_j]  = 
\Ex\Ey\left[\gD_i\cdot\gD_j \ind_{\{Y_r=X_r\}} \right] \Pxy(Y_r=X_r)^{-1}}\nonumber\\
 & = & \Ex \Big[ \Ey\left[\gD_i\cdot\gD_j \left|Y_r=X_r\right.\right]  \Py(Y_r=X_r) \Big] \Pxy(Y_r=X_r)^{-1} \nonumber\\
 & = & \frac{1}{r-1}\left( \Ex\left[\frac{\norm{X_{r}}^2}{r}\Py(Y_r=X_r)\right] \Pxy(Y_r=X_r)^{-1} \right.\nonumber\\
   & & \hspace{2cm} \left. - \Ex\Ey\left[ \norm{\gD_1}^2 \ind_{\{Y_{r}=X_r\}} \right]\Pxy(Y_r=X_r)^{-1} \right) \nonumber\\
 & = & \frac{1}{r-1}\left( \Ex\left[\frac{\norm{X_{r}}^2}{r} \Py(Y_r=X_r)\right] \Pxy(Y_r=X_r)^{-1} - \Ex\Ey\left[ \norm{\gD_1}^2 \left|Y_r=X_r\right. \right] \right).\nonumber
\end{eqnarray}

Next, we study the 
 asymptotic behavior of $A(r)$ and we prove \eqref{A(r)} with 
$C_{X,Y}=tr (\gS_Y)- tr \left( (\gS_X^{-1}+ \gS_Y^{-1})^{-1}\right)$.
Note that $tr(\gS_Y)=\bbE^Y (||Y_1||^2):=\gs_Y^2$.
The fact that 
$C_{X,Y}>0$ is just a consequence of the fact that, if $A$ and $B$ 
are two positive-definite matrices, one has that $A-B$ is positive definite 
if and only if 
$B^{-1}- A^{-1}$ is \cite[Cor. 7.7.4(a)]{hornjohnson}.

To prove \eqref{A(r)}, it is enough to show that
\begin{equation}
  \Ex\Ey\left[ \norm{\gD_1}^2 \left|Y_r=X_r\right. \right] 
\stackrel{r\to\infty}{\rightarrow} \Ex\Ey\left[ \norm{\gD_1}^2 \right] =\gs_Y^2,
\label{asympA part1}
\end{equation}
and that
\begin{equation} 
  B(r):=\frac{\Ex\left[\frac{\norm{X_r}^2}{r} \Py(Y_r=X_r)\right]}{\Pxy(X_r=Y_r)} \stackrel{r\to\infty}{\rightarrow} tr 
\left( (\gS_X^{-1}+ \gS_Y^{-1})^{-1}\right).
\label{asympA part2}
\end{equation}


To prove \eqref{asympA part1}, write
\begin{eqnarray}
  \Ex\Ey\left[ \norm{\gD_1}^2 \left|Y_r=X_r\right. \right] 
  = \Ey\left[ \norm{\gD_1}^2 \Px(X_r=Y_r) \right]\Pxy(X_r=Y_r)^{-1} \label{Exy gDi xr=yr}.
\end{eqnarray}
We know from Assumption \ref{hypwright} and from simple large
deviation bounds that $\Py(\norm{\gD_1}> r^{1/4})$ and $\Py\left(\norm{Y_r}
> \frac{r^{3/5}}{2}\right)$ decay faster than any inverse power of $r$ for
$r\to\infty$, so that 
\begin{eqnarray}
 \Ex\Ey\left[ \norm{\gD_1}^2 \Px(X_r=Y_r) \right]  
 =(1+o(1))\Ex\Ey\left[ \norm{\gD_1}^2 \Px(X_r=Y_r)
        \ind_{\{\norm{\gD_1}\leq  r^{1/4}\}} \ind_{\{\norm{Y_r} \leq \frac{r^{3/5}}{2}\}} \right] \nonumber
\end{eqnarray}
We now decompose the expectation according to the values of $Y_r$ and $\gD_1$
and we use the Local Limit Theorem, Proposition \ref{prop:localCLT}
(observe that $||Y_r-\Delta_1||\le r^{3/5}$):
\begin{eqnarray}
 \lefteqn{\sum_{\norm{x}\leq (1/2)r^{3/5}}
\sum_{\|x_1\|\leq r^{1/4}} \norm{x_1}^2 \Py(Y_1=x_1)\Py(Y_{r-1}=x-x_1)\Px(X_r=x) } \nonumber\\
  & = & (1+o(1))\sum_{\|x_1\|\leq r^{1/4}} 
\norm{x_1}^2 \Py(Y_1=x_1)
\sum_{\norm{x}\leq \frac{r^{3/5}}{2}}\frac{c_X c_Y}{r^d}
e^{-\frac{1}{2r} (x-x_1)\cdot\left(\gS_Y^{-1}(x-x_1)\right)} 
 e^{ -\frac{1}{2r} x\cdot \left(\gS_X^{-1}x\right)} 
\nonumber\\
  & = & (1+o(1))\gs_Y^2 \left( \sum_{\norm{x}\leq r^{3/5}}\frac{c_X c_Y}{r^{d}} e^{-\frac{1}{2r} x\cdot\left((\gS_X^{-1}+\gS_Y^{-1})x\right)} \right)\label{Exy gDi sum},
\end{eqnarray}
where $c_X=(2\pi)^{-d/2}(\det \Sigma_X)^{-1/2}$ and 
similarly for $c_Y$ (the constants are different in the case of simple random walks: see Remark~\ref{CLTsimple}), and where we used the fact 
$(x-x_1)\cdot
\left(\gS_Y^{-1}(x-x_1)\right)=x\cdot\left(\gS_Y^{-1}x\right)+o(r)$
uniformly for the values of $x,x_1$ we are considering.

Using the same reasoning, 
we also have (with the same constants $c_X$ and $c_Y$)
\begin{eqnarray}
\Pxy(X_r=Y_r)  & = &  (1+o(1))\sum_{\norm{x}\leq r^{3/5}}\Py(Y_r=x)\Px(X_r=x) \nonumber\\
  & = & (1+o(1)) \sum_{\norm{x}\leq r^{3/5}} \frac{c_X c_Y}{ r^{d}} e^{-\frac{1}{2r} x\cdot\left((\gS_X^{-1}+\gS_Y^{-1})x\right)}.
\label{Pxr=yr}
\end{eqnarray}
Putting this together with (\ref{Exy gDi xr=yr}) and (\ref{Exy gDi sum}), we have
the asymptotic behavior of the term in~(\ref{asympA part1}).

\begin{rem}\rm
\label{rem:asympAr}
For the
purposes of Section~\ref{proof variance}, we remark that with the
same method one can prove that, for any polynomial
$Q$ of $\norm{\gD_1},\norm{\gD_2},\norm{\gD_3}$, one has
\begin{equation}
 \Ex\Ey\big[ Q(\norm{\gD_1}, \norm{\gD_2}, \norm{\gD_3}) \left| Y_r=X_r \right. \big] \stackrel{r\to\infty}{\rightarrow} 
    \Ey\big[ Q(\norm{\gD_1}, \norm{\gD_2}, \norm{\gD_3}) \big].
\label{asymp Ect Q(gD1)}
\end{equation}
\end{rem}

To deal with the term $B(r)$ in~(\ref{asympA part2}), we use the same
method: we know that $\Px(||X_r|| > r^{3/5})$ decays faster than any
inverse power of $r$, so decomposing according to the values of $X_r$
and applying the Local Limit Theorem we have
\begin{equation}
\Ex\left[\frac{\norm{X_r}^2}{r} \Py(Y_r=X_r)\right] 
= 
(1+o(1))  \frac{ c_Y c_X}{r^d} 
  \sum_{\norm{x}\leq r^{3/5}}\frac{\norm{x}^2}{r} e^{ -\frac{1}{2r} x\cdot \left(\gS_Y^{-1}x\right)} e^{ -\frac{1}{2r} x\cdot\left(\gS_X^{-1}x\right)}.
\end{equation}
Together with \eqref{Pxr=yr}, we finally get 
\begin{equation}
B(r) = (1+o(1)) \frac{\sum_{\norm{x}\leq r^{3/5}}\frac{\norm{x}^2}{r}
e^{ -\frac{1}{2r} x\cdot\left((\gS_Y^{-1}+\gS_X^{-1})x\right)} }{ \sum_{\norm{x}\leq r^{3/5}}   e^{ -\frac{1}{2r} x\cdot\left((\gS_Y^{-1}+\gS_X^{-1})x\right)} }=(1+o(1))E\left[\norm{ \cN }^2\right],
\label{expressAn}
\end{equation}
where $\cN\sim \cN\left(0, (\gS_Y^{-1}+\gS_X^{-1})^{-1}\right)$ is a centered Gaussian vector of covariance matrix 
${(\gS_Y^{-1}+\gS_X^{-1})^{-1}} $. Therefore, $E\left[\norm{ \cN }^2\right]=tr \left( (\gS_Y^{-1}+\gS_X^{-1})^{-1}\right)$ and \eqref{asympA part2} is proven.

\qed

\medskip
Given Lemma \ref{th:Ar}, we can resume the proof of Lemma \ref{EctF big},
and lower bound the average 
 $\Ect[F_1]$. 
Recalling \eqref{F_1} and the fact that we reduced to the case $a=0$, 
we get
\begin{eqnarray}
\Ect[F_1] & = & \sum_{n=1}^{R_b} \left(\sum_{\tau_{n-1}<i,j\le \tau_n}M_{ij}\right) A(\gdt_n),
\label{bound1 EctF}
\end{eqnarray}
where $\gdt_n:=\tau_n-\tau_{n-1}$.
Using the definition~(\ref{defM})
of $M$, we see that there exists a constant $c>0$ such
that for $1<m\le L$
\begin{equation}
\sum_{i,j=1}^m M_{ij}  \ge \frac{c}{\sqrt{L\log L}} m^{3/2}.
\end{equation}
On the other hand, thanks to Lemma \ref{th:Ar}, there exists some $r_0>0$ and two constants $c$ and
$c'$ such that $A(r)\geq
\frac{c}{r}$ for $r\geq r_0$, and $A(r)\geq -c'$ for every
$r$.  Plugging this into~(\ref{bound1 EctF}), one gets
\begin{equation}
\sqrt{L\log L}\,\Ect[F_1]  \geq  c \sum_{n=1}^{R_b} \sqrt{\gdt_n}\ind_{\{\gdt_n \hspace{-0.05cm}\geq \hspace{-0.05cm}r_0\}}
         -c' \sum_{n=1}^{R_b}(\gdt_n)^{3/2}\ind_{\{\gdt_n \hspace{-0.05cm}\leq \hspace{-0.05cm}r_0\}} 
   \geq  c \sum_{n=1}^{R_b} \sqrt{\gdt_n} - c'R_b. \label{bound2 EctF}
\end{equation}
Therefore, we get for any positive $B>0$ (independent of $L$)
\begin{eqnarray}
\lefteqn{\PL\left(\Ect[F_1] \leq g\sqrt{\log L}\right)  \leq  \PL \left(\frac{1}{\sqrt{L\log L}}\left(c\sum_{n=1}^{R_b} \sqrt{\gdt_n}
        - c'  \,R_b\right) \leq g\sqrt{\log L}\right) } 
\nonumber\\
   &\hspace{2cm} \leq & \PL\left(\frac{1}{\sqrt{L\log L}}\left( c\sum_{n=1}^{R_b} \sqrt{\gdt_n}
         - c'\sqrt{L}B\right)\leq g \sqrt{\log L}\right) +\PL\left(R_b> B\sqrt{L} \right) \nonumber\\
   &\hspace{2cm} \leq & \PL \left( \sum_{n=1}^{R_{b/2}} \sqrt{\gdt_n}\le (1 + o(1))\frac1c g\sqrt{L}\log L \right) + \PL(R_b > B \sqrt{L}). \label{bound3 EctF}
\end{eqnarray}
Now we show that for $B$ large enough, and $L\geq L_0(B)$,
\begin{equation}
\PL(R_b > B \sqrt{L})\leq \gz/2,
\label{PtauWbig}
\end{equation}
where $\gz$ is the constant which appears in the statement of Lemma \ref{EctF big}.
We start with getting rid of the conditioning in $\PL$ (recall $\PL(\cdot) = \PL(\cdot|b\in\tau)$ since we reduced to the case $a=0$).
If $R_{b}>B \sqrt{L}$,
then either $\left|\tau \cap\{1,\ldots, b/2\}\right| $ or
$\left|\tau \cap\{b/2+1,\ldots, b\}\right|$ exceeds $ \frac{B}{2} \sqrt{L}$.
Since both random variables have the same law under $\PL$,
 we have
\begin{equation}
\PL(R_b> B \sqrt{L}) \leq 2\PL\left(R_{b/2}>\frac{B}{2} \sqrt{L}\right)\le 2 c\P\left(R_{b/2}>\frac{B}{2} \sqrt{L}\right),
\end{equation}
where in the second inequality we applied 
Lemma~\ref{decond2}.
Now, we can use the Lemma~\ref{conv in law W} in the Appendix, to get that (recall $b\le L$)
\begin{equation}
\P\left(R_{b/2}>\frac{B}{2} \sqrt{L}\right) \le \P\left(R_{L/2}>\frac{B}{2} \sqrt{L}\right) 
\stackrel{L\to\infty}{\rightarrow}
  \P\left(\frac{\left| \cZ \right|}{\sqrt{2\pi}}\geq B\frac{ c_K}{\sqrt{2}} \right),
\end{equation}
with $\cZ$ a standard Gaussian random variable and $c_K$ the constant such that $K(n)\sim c_K n^{-3/2}$. The inequality
\eqref{PtauWbig} then follows 
for $B$ sufficiently large, and $L\geq L_0(B)$.

We are left to prove that for $L$ large enough and $g$ small enough
\begin{equation}
\PL \left( \sum_{n=1}^{R_{b/2}} \sqrt{\gdt_n} \leq  \frac gc\sqrt{L}\log L \right) \leq \gz/2.
\end{equation}
The conditioning in $\PL$ can be eliminated again via Lemma~\ref{decond2}. Next, one notes that for any given $A>0$ (independent of $L$)
\begin{equation}
\P \left( \sum_{n=1}^{R_{b/2}} \sqrt{\gdt_n}\leq  \frac gc\sqrt{L}\log L \right) 
   \leq  \P \left( \sum_{n=1}^{A\sqrt{L}} \sqrt{\gdt_n}\leq \frac gc\sqrt{L}\log L \right) + \P\left(R_{b/2}< A\sqrt{L}\right).
\label{eq:sum gdt}
\end{equation}
Thanks to the Lemma~\ref{conv in law W} in Appendix and to $b\ge \gep L$, we have
$$\limsup_{L\to\infty}\P\left(\frac{R_{b/2}}{\sqrt{L}}<A \right)\le\P\left(\frac{\left| \cZ \right|}{\sqrt{2\pi}}< A c_K\sqrt{\frac2\gep}\right),$$
which can be arbitrarily small if $A=A(\gep)$ is small enough, for $L$ large.
We now deal with the other term in \eqref{eq:sum gdt}, using the exponential Bienaym\'e-Chebyshev inequality
(and the fact that the $\gdt_n$ are i.i.d.):
\begin{eqnarray}
\P \left( \frac{1}{\sqrt{L\log L}}\sum_{n=1}^{A\sqrt{L}} \sqrt{\gdt_n} < \frac gc\sqrt{\log L} \right) 
  \leq  e^{(g/c)\sqrt{\log L}} \E \left[\exp\left( - \sqrt{ \frac{\tau_1}{L\log L} } \right)\right]^{A\sqrt{L}}. \label{bienayme expo}
\end{eqnarray}
To estimate this expression, we remark that, for $L$ large enough,
\begin{eqnarray}
\E \left[1-\exp\left( - \sqrt{ \frac{\tau_1}{L\log L}  } \right)\right] & = & 
\sum_{n=1}^{\infty} K(n)\left(1-e^{- \sqrt{ \frac{n}{L\log L} } }\right) \nonumber\\
 & \geq  & c'\sum_{n=1}^{\infty} \frac{1-e^{- \sqrt{ \frac{n}{L\log L}} }}{n^{3/2}}\ge c''\sqrt{\frac{\log L}L},
\end{eqnarray}
where the last inequality follows from keeping only the terms with $n\le L$ in the sum, and noting that in this range 
$1-e^{- \sqrt{ \frac{n}{L\log L}}}\ge c\sqrt{n/(L\log L)}$.
Therefore,
\begin{eqnarray}
\E \left[\exp\left( - \sqrt{ \frac{\tau_1}{L\log L} } \right)\right]^{A \sqrt{L}}   \leq  \left(1-c''\sqrt{\frac{\log L}{L}} \right)^{A \sqrt{L}} 
 \le  e^{-c''A\sqrt{\log L}},
\end{eqnarray}
and, plugging this bound in the inequality~(\ref{bienayme expo}), we get
\begin{equation}
\bP \left( \frac{1}{\sqrt{L\log L}}\sum_{n=1}^{A\sqrt{L}} \sqrt{\gdt_n} \leq 
\frac gc\sqrt{\log L}\right) \leq e^{[(g/c)-c''A]\sqrt{\log L}},
\end{equation}
that goes to $0$ if $L\to\infty$, provided that $g$ is small enough. This concludes the proof of Lemma~\ref{EctF big}.
\qed

\subsection{Proof of the Lemma~\ref{variance} }
\label{proof variance}
We can write
\begin{equation}
F_1-\Ect[F_1] = S_1+S_2:=\sum_{i\ne j=a+1}^{b} M_{ij} D_{ij}+
\sum_{i\ne j} ' M_{ij} D_{ij}
\end{equation}
where we denoted
\begin{equation}
D_{ij}=\gD_i \cdot\gD_j-\Ect[\gD_i\cdot \gD_j]
\end{equation}
and $\stackrel{\prime}{\sum}$ stands for the sum over all $1\le i\ne j\le L$ such that
either $i$ or $j$ (or both) do not fall into $\{a+1,\ldots,b\}$.
This way, we have to estimate
\begin{eqnarray}
\label{sum variance}
\lefteqn{\Ect[(F_1-\Ect[ F_1])^2] \le
2\Ect[S_1^2]+2\Ect[S_2^2]}\\\nonumber
&\hspace{2cm}=& 
2\sum_{i\ne j=a+1}^b\sum_{k\ne l=a+1}^{b} M_{ij}M_{kl} \Ect[D_{ij}D_{kl}]+2
\sum'_{i\ne j}\sum'_{k\ne l}M_{ij}M_{kl} \Ect[D_{ij}D_{kl}].
\end{eqnarray}

\begin{rem}\rm
\label{rem:ijkl}
We easily deal with the part of the sum where $\{i,j\}=\{k,l\}$. In fact,
we trivially bound
  $\Ect\left[(\gD_i\cdot\gD_j)^2\right]\le\Ect\left[\norm{\gD_i}^4\right]^{1/2}
    \left[\norm{\gD_j}^4 \right]^{1/2}$.
Suppose for instance that $\tau_{n-1}<i\le \tau_n$ for some
$R_a<n\le R_b$: in this case $\Ect\left[\norm{\gD_i}^4 \right]$ converges
for $\tau_n-\tau_{n-1}\to\infty$  to $\Ey[\norm{\gD_1}^4]$  thanks to 
  (\ref{asymp Ect Q(gD1)}). If, on the other hand, $i\notin \{a+1,\ldots,b\}$, we know
that $\Ect\left[\norm{\gD_i}^4\right]$ equals exactly  $\Ey\left[\norm{\gD_1}^4\right]$.

As a consequence, we have the following inequality, valid for every $1\le i,j \le L$:
  \begin{equation}
  \Ect\left[(\gD_i\cdot\gD_j)^2\right]\le c.
  \end{equation}
And then
\begin{equation}
\sum_{i\ne j=1}^L \sum_{\{k,l\}=\{i,j\}} M_{ij}M_{kl} \Ect[D_{ij}D_{kl}]\leq c \sum_{i\ne j=1}^L M_{ij}^2 \leq c'
\end{equation}
since the Hilbert-Schmidt norm of $M$ was chosen to be finite.
\end{rem}

{\bf Upper bound on $\Ect[S_2^2]$}.
This is the easy part, and this term will be shown to be bounded even
without taking the average over $\PL$.\\
We have to compute $\stackrel{\prime}{\sum}_{i\ne j}\stackrel{\prime}{\sum}_{k\ne l}M_{ij}M_{kl}
\Ect[D_{ij}D_{kl}]$.  Again, thanks to \eqref{dif intervals}-\eqref{ident}, 
we have $\Ect[D_{ij}D_{kl}]\neq 0$ only in the following
case (recall that thanks to Remark~\ref{rem:ijkl} we can disregard the
case $\{i,j\}=\{k,l\}$):
\begin{eqnarray}
  i=k\notin\{a+1,\ldots,b\} \mbox{\;and\;} \tau_{n-1}<j\ne l\le \tau_n \mbox{\;for some
\;}
R_a< n\le R_b.
\label{eq:casoikl}
\end{eqnarray}
One should also consider the cases where $i$ is interchanged with $j$
and/or $k$ with $l$.
Since we are not following constants, we do not keep track of the associated combinatorial
factors.
%
%
Under the assumption \eqref{eq:casoikl}, $\Ect[\gD_i\cdot\gD_j]=\Ect[\gD_i\cdot\gD_l]=0$ (cf. \eqref{dif
intervals}) and we will show that
\begin{eqnarray}
  \label{eq:strazio}
\Ect[D_{ij}D_{il}]= \Ect[(\gD_i\cdot\gD_j)(\gD_i\cdot\gD_l)]  \le \frac cr
\end{eqnarray}
where $r=\tau_n-\tau_{n-1}=\gdt_n$. Indeed,
using \eqref{dif intervals}-\eqref{ident}, we get
\begin{eqnarray}
  \nonumber
  \Ect[(\gD_i\cdot\gD_j)(\gD_i\cdot\gD_l)] &=&
  \sum_{\nu,\mu=1}^3\bbE^Y[\Delta_i^{(\nu)}\Delta_i^{(\mu)}]
  \bbE^X\bbE^Y[\Delta_{j-\tau_{n-1}}^{(\nu)}\Delta_{l-\tau_{n-1}}^{(\mu)}|X_{\tau_n-\tau_{n-1}}=Y_{\tau_n-\tau_{n-1}}]\\
  &=&\sum_{\nu,\mu=1}^3 \Sigma_Y^{\nu\mu}\,
  \bbE^X\bbE^Y\left[\Delta_{j-\tau_{n-1}}^{(\nu)}\Delta_{l-\tau_{n-1}}^{(\mu)}|X_r=Y_r\right].
\end{eqnarray}
In the remaining expectation, we can assume without loss of generality that $\tau_{n-1}=0,\tau_n=r$. 
Like for instance in the proof of \eqref{A(r)}, one writes
\begin{eqnarray}
\label{454}
  \bbE^X\bbE^Y\left[\Delta_j^{(\nu)}\Delta_l^{(\mu)}|X_{r}=Y_{r}\right]=\frac{\bbE^X \left[\bbE^Y\big[ \Delta_j^{(\nu)}\Delta_l^{(\mu)}|Y_r=X_r\big] \bbP^Y(Y_r=X_r)\right]
}{\bbP^{X-Y}(X_r=Y_r)}
\end{eqnarray}
and
\begin{eqnarray}
\label{455}
  \bbE^Y\left[\left.\Delta_j^{(\nu)}\Delta_l^{(\mu)}\right|Y_r=X_r\right] = \frac1{r(r-1)}
X_r^{(\nu)}X_r^{(\mu)}
-
  \frac1{r-1}\bbE^Y[ \Delta_j^{(\nu)}\Delta_j^{(\mu)}|Y_r=X_r].
\end{eqnarray}
An application of the Local Limit Theorem like in \eqref{asympA part1},
\eqref{asympA part2} then leads to 
\eqref{eq:strazio}.

We are now able to bound
\begin{eqnarray}
\label{eq:S22}
\lefteqn{\bbE_\tau\left[S_2^2\right]=c \sum_{i\notin\{a+1,\ldots,b\}} \sum_{n=R_a+1}^{R_b} \sum_{\tau_{n-1}< j,l \leq \tau_n } M_{ij}M_{il} \Ect[D_{ij}D_{il}]}\nonumber\\
&\hspace{2cm}\leq & \frac{c}{L\log L}\sum_{i\notin\{a+1,\ldots,b\}} \sum_{n=R_a+1}^{R_b}
\sum_{\tau_{n-1}< j,l \leq \tau_n } \frac{1}{\sqrt{|i-j|}} \frac{1}{\sqrt{|i-l|}} \frac{1}{\gdt_n }.
\end{eqnarray}
Assume for instance that $i>b$ (the case $i\le a$ can be treated similarly): 
\begin{eqnarray}
\nonumber
\lefteqn{\frac{c}{L\log L}\sum_{i> b} \sum_{n=R_a+1}^{R_b}
\sum_{\tau_{n-1}< j,l \leq \tau_n } \frac{1}{\sqrt{i-j}} \frac{1}{\sqrt{i-l}} \frac{1}{\gdt_n }   } \\\nonumber
&\le& \frac{c}{L\log L}\sum_{i> b} \sum_{n=R_a+1}^{R_b}\sum_{\tau_{n-1}< j,l \leq \tau_n } \frac{1}{(i-\tau_n)\gdt_n } 
\le \frac{c}{L\log L} (b-a)\sum_{i=1}^L\frac1i\le c'.
\end{eqnarray}

{\bf Upper bound on $\Ect[S_1^2]$}.
Thanks to time translation invariance,
one can reduce to the case $a=0$.
We have to distinguish various cases (recall
Remark~\ref{rem:ijkl}: we assume that $\{i,j\}\neq\{k,l\}$).

\begin{enumerate}
\item Assume that
  $\tau_{n-1}<i,j\le \tau_{n} $, $\tau_{m-1}<k,l\le \tau_m$, 
with $1\le n\neq m\le R_b$.
  Then, thanks to \eqref{dif intervals}, we get
  $\Ect[D_{ij}D_{kl}]=\Ect[D_{ij}]\Ect[D_{kl}]=0$, because
  $\Ect[D_{ij}]=0$.
For similar reasons, one has that  $\Ect[D_{ij}D_{kl}]=0$ if one of the indexes, say $i$, belongs to 
one of the intervals $\{\tau_{n-1}+1,\ldots,\tau_n\}$, and the other three do not.

\item Assume that $\tau_{n-1}<i,j,k,l\le\tau_n$ for some $n\le R_b$.
  Using \eqref{ident}, we have
$$\Ect[D_{ij} D_{kl}] = 
\Ey\Ex\left[D_{ij} D_{kl}\left|X_{\tau_{n-1}}=Y_{\tau_{n-1}}, X_{\tau_{n}}=Y_{\tau_{n}} \right.\right],$$
and with a time translation we can reduce to the case $n=1$ (we call $\tau_1=r$).
Thanks to the computation of
$\Ect[\gD_i\cdot\gD_j]$ in  Section~\ref{lemma EctF big}, we see that
$\Ect[\gD_i\cdot\gD_j]=\Ect[\gD_k\cdot\gD_l]=-A(r)$ so that
\begin{equation}
\Ect[D_{ij}D_{kl}] = \Ect[(\gD_i\cdot\gD_j)(\gD_k\cdot\gD_l)]-A(r)^2 \leq  \Ect[(\gD_i\cdot\gD_j)(\gD_k\cdot\gD_l)].
\end{equation}
\begin{enumerate}


  \item If $i=k$, $j\neq l$ (and $\tau_{n-1}<i,j,l\le\tau_n$ for some $n\le R_b$), then 
    \begin{eqnarray}
           \label{eq:ijl}
\bbE_\tau[(\Delta_i\cdot \Delta_j)(\Delta_i\cdot\Delta_l)]\le \frac c{\Delta\tau_n}.
    \end{eqnarray}
The computations are similar to those we did in Section~\ref{lemma EctF big} for the computation of $\Ect[\gD_i\cdot\gD_j]$. See Appendix \ref{app:ijl} for details.

 \item If $\{i,j\}\cap \{k,l\}=\emptyset$ (and $\tau_{n-1}<i,j,k,l\le\tau_n$ for some $n\le R_b$), one gets
    \begin{eqnarray}
           \label{eq:ijkl}
\bbE_\tau[(\Delta_i\cdot \Delta_j)(\Delta_k\cdot\Delta_l)]\le \frac c{(\Delta\tau_n)^2}.
    \end{eqnarray}
See  Appendix \ref{app:ijkl} for a (sketch of) the proof, which is analogous to that of \eqref{eq:ijl}.

\end{enumerate}

\item The only remaining case is that where
  $i\in\{\tau_{n-1}+1,\ldots,\tau_{n}\}$,
  $j\in\{\tau_{m-1}+1,\ldots,\tau_{m}\}$ with $m\neq n\le R_b$, and
  each of these two intervals contain two indexes in $i,j,k,l$.  Let us suppose for
  definiteness $n<m$ and $k\in \{\tau_{n-1}+1,\ldots,\tau_{n}\}$.
  Then $\Ect[\gD_i\cdot\gD_j]=\Ect[\gD_k\cdot\gD_l]=0$ (cf. Lemma
  \ref{th:Ar}), and $\Ect[D_{ij}D_{kl}] =  \Ect[(\gD_i\cdot\gD_j)(\gD_k\cdot\gD_l)]$. 
We will prove in Appendix \ref{app:ijklmn} that 
\begin{equation}
\Ect[(\gD_i\cdot\gD_j)(\gD_k\cdot\gD_l)]  \leq   \frac c{\gdt_n \gdt_m}.
\label{ijklmn}
\end{equation}
\end{enumerate}

We are now able to compute $\bbE_\tau[S_1^2]$. We consider first 
the contribution of the terms whose indexes
$i,j,k,l$ are all in the same interval $\{\tau_{n-1}+ 1,\ldots,\tau_{n}\}$, i.e. case (2) above. Recall that we drop the
terms $\{i,j\}=\{k,l\}$ (see Remark~\ref{rem:ijkl}):
\begin{eqnarray}
  \lefteqn{\sumtwo{\tau_{n-1}<i,j,k,l\le \tau_n }{\{i,j\}\neq\{k,l\}} M_{ij}M_{kl} \Ect[D_{ij}D_{kl}] \leq
     \frac{c}{\gdt_n} \sumtwo{l\in\{i,j\} \ or\ k\in\{i,j\}  }
    {\tau_{n-1}<i,j,k,l\le \tau_n } M_{ij}M_{kl} + \frac{c}{\gdt_n^2} \sumtwo{\{i,j\}\cap\{k,l\}=\emptyset }
    {\tau_{n-1}<i,j,k,l\le \tau_n }  M_{ij}M_{kl}} \nonumber\\
  & \hspace{2.5cm}\leq & \frac {c'}{L\log L} \left[ \frac{1}{\gdt_n} \sum_{1\le i<j<k\le\gdt_n} \frac{1}{\sqrt{j-i}}\frac{1}{\sqrt{k-j}} + 
\frac{1}{\gdt_n^2} \left(\sum_{1\le i<j\le \gdt_n} \frac{1}{\sqrt{j-i}}\right)^2 \right] \nonumber\\
  & \hspace{2.5cm}\leq & \frac{c''}{L\log L} \gdt_n.
\end{eqnarray}

Altogether, we see that
\begin{eqnarray}
\lefteqn{\sum_{i\ne j=1}^b\sumtwo{k\ne l=1}{\{i,j\}\neq\{k,l\}}^{b} M_{ij}M_{kl} \Ect[D_{ij}D_{kl}]\ind_{\{\exists n\le R_b:
i,j\in \{\tau_{n-1}+1,\ldots,\tau_n\}\}} } \nonumber\\
 &\hspace{2.5cm}= & \sum_{n=1}^{R_b} \sumtwo{\tau_{n-1}<i,j,k,l\le\tau_n}{\{i,j\}\neq\{k,l\}} M_{ij}M_{kl} \Ect[D_{ij}D_{kl}]  \leq
 \frac{c}{L \log L} \sum_{n=1}^{R_b}  \gdt_n \le \frac c{\log L}. \label{case i,j same interval}
\end{eqnarray}

Finally, we consider the contribution to $\bbE_\tau[S_1^2]$ coming from the terms of point (3). We have (recall that
$n<m$)
\begin{eqnarray}
\sumtwo{\tau_{n-1}<i,k\le\tau_n }{\tau_{m-1}<j,l\le \tau_m} M_{ij}M_{kl} \Ect[D_{ij}D_{kl}]
  \leq  \frac{ c}{L\log L}\frac{1}{\gdt_n\gdt_m} \sumtwo
{\tau_{n-1}<i,k\le\tau_n }{\tau_{m-1}<j,l\le \tau_m} \frac{1}{\sqrt{j-i}}\frac{1}{\sqrt{l-k}}.
\end{eqnarray}
But as $j> \tau_{m-1}$
\begin{equation}
\sum_{\tau_{n-1}< i\leq \tau_{n}}\frac{1}{\sqrt{j-i}} \leq \sum_{\tau_{n-1}< i\leq \tau_{n}}\frac{1}{\sqrt{\tau_{m-1}-i+1}} 
\leq c \left(\sqrt{\tau_{m-1}-\tau_{n-1}}-\sqrt{\tau_{m-1}-\tau_{n}}\right),
\end{equation}
and as $k\leq \tau_{n}$
\begin{equation}
\sum_{\tau_{m-1}< l\leq \tau_{m}}\frac{1}{\sqrt{l-k}} \leq \sum_{\tau_{m-1}< l\leq \tau_{m}}\frac{1}{\sqrt{l-\tau_{n}}} \leq c \left(\sqrt{\tau_{m}-\tau_{n}}-\sqrt{\tau_{m-1}-\tau_{n}}\right),
\end{equation}
so that
\begin{equation}
\sumtwo{\tau_{n-1}<i,k\le\tau_n }{\tau_{m-1}<j,l\le \tau_m} M_{ij}M_{kl} \Ect[D_{ij}D_{kl}] \leq  \frac{c}{L \log L} \left(\sqrt{T_{nm}+\gdt_n}-\sqrt{T_{nm}}\right)\left(\sqrt{T_{nm}+\gdt_m}-\sqrt{T_{nm}}\right),
\end{equation}
where we noted
$T_{nm} =\tau_{m-1}-\tau_{n}.$
Recalling \eqref{case i,j same interval} and the definition 
 \eqref{sum variance} of 
$S_1$,
 we can finally write
\begin{eqnarray}
\lefteqn{\EL \left[ \bbE_\tau[S_1^2]\right] 
\le 
c\left(
1+  \EL\left[ \sum_{n=1}^{R_b-1} \sum_{n<m\le R_b}\sumtwo{\tau_{n-1}<i,k\le\tau_n }{\tau_{m-1}<j,l\le \tau_m} M_{ij}M_{kl} \Ect[D_{ij}D_{kl}] \right]\right) }\nonumber\\
 & \leq &c+ \frac{c}{L\log L} \EL\left[\sum_{1\le n<m\le R_b} \left(\sqrt{T_{nm}+\gdt_n}-\sqrt{T_{nm}}\right)\left(\sqrt{T_{nm}+\gdt_m}-\sqrt{T_{nm}}\right) \right]. \nonumber
\end{eqnarray}
The remaining average can be estimated via the following Lemma.
\begin{lemma}
\label{th:derniert}
There exists a constant $c>0$ depending only on $K(\cdot)$, such that
\begin{equation}
\EL\left[\sum_{1\le n<m\le R_b}\left(\sqrt{T_{nm}+\gdt_n}-\sqrt{T_{nm}}\right)\left(\sqrt{T_{nm}+\gdt_m}-\sqrt{T_{nm}}\right)\right] \leq c L(\log L)^{7/4}.
\label{dernier terme}
\end{equation}
\end{lemma}
Of course this  implies that $\EL\bbE_\tau[S_1^2]\le c (\log L)^{3/4}$,
which together with \eqref{eq:S22} implies the claim of Lemma \ref{variance}.
\qed

\medskip

{\sl Proof of Lemma \ref{th:derniert}}.
One has the inequality
\begin{equation}
  \left(\sqrt{T_{nm}+\gdt_n}-\sqrt{T_{nm}}\right)\left(\sqrt{T_{nm}+\gdt_m}-\sqrt{T_{nm}}\right)\leq  \sqrt{\gdt_n} \sqrt{\gdt_m},
\end{equation}
which is a good approximation when $T_{nm}$ is not that large compared
with $\gdt_n$ and $\gdt_m$, and 
\begin{equation}
\left(\sqrt{T_{nm}+\gdt_n}-\sqrt{T_{nm}}\right)\left(\sqrt{T_{nm}+\gdt_m}-\sqrt{T_{nm}}\right)\leq c \frac{\gdt_n \gdt_m}{T_{nm}} ,
\end{equation}
which is accurate when $T_{nm}$ is large.
We use these bounds to cut the expectation~(\ref{dernier
  terme}) into two parts, a term where $m-n\leq H_L$ and one where $m-n>
H_L$, with $H_L$ to be chosen later:
\begin{eqnarray}
\lefteqn{\EL\left[\sum_{n=1}^{R_b}\sum_{m=n+1}^{R_b}\left(\sqrt{T_{nm}+\gdt_n}-\sqrt{T_{nm}}\right)\left(\sqrt{T_{nm}+\gdt_m}-\sqrt{T_{nm}}\right)\right]} \nonumber\\
 & \hspace{1cm}\leq & \EL\left[\sum_{n=1}^{R_b}\sum_{m=n+1}^{(n+H_L)\wedge R_b} 
\sqrt{\gdt_n}\sqrt{\gdt_{m}} \right] + c\, \EL\left[\sum_{n=1}^{R_b}\sum_{m = n+H_L+1}^{R_b} \frac{\gdt_n \gdt_m}{T_{nm}}\right]. \label{dernier terme part1}
\end{eqnarray}

We claim  that there exists a constant $c$ such that for every
 $l\geq 1$,
\begin{equation}
\label{112}
\EL\left[\sum_{n=1}^{R_b-l}\sqrt{\gdt_n}\sqrt{\gdt_{n+l}} \right] \leq c \sqrt{L}(\log L)^{2+\frac{1}{12}}
\end{equation}
(the proof is given later).
Then the first term in the right-hand side of (\ref{dernier terme part1}) is
\begin{equation}
  \EL\left[\sum_{n=1}^{R_b}\sum_{m=n+1}^{(n+H_L)\wedge R_b} \sqrt{\gdt_n}\sqrt{\gdt_{m}} \right]  =  \sum_{l=1}^{H_L} \EL\left[\sum_{n=1}^{R_b-l}\sqrt{\gdt_n}\sqrt{\gdt_{n+l}}\right] 
   \leq  c H_L \sqrt{L}(\log L)^{2+1/12}. \nonumber
\end{equation}
If we choose
$H_L= 
\sqrt{L}(\log L)^{-1/3} ,$
we get  from ~(\ref{dernier terme part1})
\begin{eqnarray}
\lefteqn{\EL\left[\sum_{n=1}^{R_b}\sum_{m=n+1}^{R_b}\left(\sqrt{T_{nm}+\gdt_n}-\sqrt{T_{nm}}\right)\left(\sqrt{T_{nm}+\gdt_m}-\sqrt{T_{nm}}\right)\right] }\label{dernier terme part2} \\
  & \hspace{2cm}\leq & c L(\log L)^{7/4} + c\, \EL\left[\sum_{n=1}^{R_b}\sum_{m = n+H_L+1}^{R_b} \frac{\gdt_n \gdt_m}{T_{nm}}\right]. \nonumber
\end{eqnarray}

As for the second term in \eqref{dernier terme part1},
recall that $T_{nm}=\tau_{m-1} - \tau_n$ and decompose the sum 
in two parts, according to whether $T_{nm}$ is larger or smaller than a
certain $K_L$ to be fixed:
\begin{eqnarray}
  \lefteqn{\EL\left[\sum_{n=1}^{R_b}\sum_{m = n+H_L+1}^{R_b} \frac{\gdt_n \gdt_m}{T_{nm}}\right]}\nonumber\\
  & \hspace{1cm}= & \EL\left[\sum_{n=1}^{R_b}\sum_{m = n+H_L+1}^{R_b} \frac{\gdt_n \gdt_m}{T_{nm}} \ind_{\{T_{nm}> K_L\}} \right] + \EL\left[\sum_{n=1}^{R_b}\sum_{m = n+H_L+1}^{R_b} \gdt_n \gdt_m \ind_{\{T_{nm}\leq K_L\}} \right] \nonumber\\
  & \hspace{1cm}\leq & \frac{1}{K_L} \EL\left[\left(\sum_{n=1}^{R_b} \gdt_n\right)^2\right] + L^2 \EL\left[\sum_{n=1}^{R_b}\sum_{m = n+H_L+1}^{R_b} \ind_{\{\tau_{n+H_L} - \tau_{n}\leq K_L\}} \right]\nonumber\\
  & \hspace{1cm}\leq & \frac{L^2}{K_L} + L^4 \PL\left( \tau_{H_L} \leq K_L  \right).
\end{eqnarray}
We now set
$K_L=L(\log L)^{-7/4}, $
so that we get in the previous inequality
\begin{equation}
\EL\left[\sum_{n=1}^{R_b}\sum_{m = n+H_L+1}^{R_b} \frac{\gdt_n \gdt_m}{T_{nm}}\right]\leq L(\log L)^{7/4} + L^4 \PL\left( \tau_{H_L} \leq K_L  \right),
\end{equation}
and we are done if we prove for instance that $\PL\left( \tau_{H_L} \leq K_L  \right) = o(L^{-4})$. Indeed,
\begin{eqnarray}
\label{WKL grand}
\PL\left( \tau_{H_L} \leq K_L \right) & = & \PL\left(R_{K_L} \geq H_L\right)\leq c \P\left(R_{K_L} \geq H_L\right)
\end{eqnarray}
where we used  Lemma~\ref{decond2} to take the conditioning off from 
$\PL:=\bP(\cdot|b\in\tau)$ (in fact, $K_L\leq b/2$ since $b\ge \gep L$).
Recalling the choices of $H_L$ and $K_L$, we get that 
$
H_L/\sqrt{K_L} =(\log L)^{13/24}$ and, combining~(\ref{WKL grand}) with  Lemma~\ref{Wn grand}, we get
\begin{equation}
\PL\left( \tau_{H_L} \leq K_L \right)  \leq c'\, e^{-c(\log L)^{13/12}} = o(L^{-4})
\end{equation}
which is what we needed.

\medskip

To conclude the proof of Lemma \ref{th:derniert}, we still have to prove \eqref{112}.
Note that
\begin{eqnarray}
  \EL\left[\sum_{n=1}^{R_b-l}\sqrt{\gdt_n}\sqrt{\gdt_{n+l}}\ind_{\{R_b>l \}}\right]   
  &= &   \EL\left[\ind_{\{R_b>l\}}\sum_{n=1}^{R_b-l}\EL\left[\sqrt{\gdt_n}\sqrt{\gdt_{n+l}}\left|R_b\right.\right]\right] \nonumber\\
  & = & \EL\left[\ind_{\{R_b>l\}}(R_b-l)
\EL\left[\sqrt{\tau_1}\sqrt{\tau_2-\tau_1}\left|R_b\right.\right]\right] \nonumber\\
  & \le & 
\EL\left[R_b\sqrt{\tau_1}\sqrt{\tau_{2}-\tau_1}\ind_{\{R_b\ge 2 \}}\right] \label{expectation first part}
\end{eqnarray}
where we used the fact that, under $\PL(\cdot|R_b=p)$ for a fixed $p$,
the law of the jumps $\{\Delta\tau_n\}_{n\le p}$ is exchangeable.
We first bound \eqref{expectation first part} when $R_b$ is large:
\begin{eqnarray}
\EL\left[R_b\sqrt{\tau_1}\sqrt{\tau_{2}-\tau_1}\ind_{\{R_b\geq \gk \sqrt{L \log L}\}}\right] & \leq & 
L^2 \PL\left(R_b\geq \gk \sqrt{L \log L} \right) \nonumber\\
& \leq & L^2 \P(b\in\tau)^{-1}  \P\left(R_b\geq \gk \sqrt{L \log L} \right).
\end{eqnarray}
In view of \eqref{eq:doney}, we have $\P(b\in\tau)^{-1}=O(\sqrt L)$.
Thanks to Lemma~\ref{Wn grand} in the Appendix, and choosing $\gk$
large enough, we get
\begin{equation}
\P\left(R_b\geq \gk \sqrt{L \log L} \right)\leq e^{-c \gk^2 \log L +o(\log L)} = o(L^{-5/2}),
\end{equation}
and therefore
\begin{equation}
\EL\left[R_b\sqrt{\tau_1}\sqrt{\tau_{2}-\tau_1}\ind_{\{R_b\geq \gk \sqrt{L \log L }\}}\right] = o(1).
\end{equation}
As a consequence, 
\begin{eqnarray}
\nonumber
\lefteqn{\EL\left[R_b\sqrt{\tau_1}\sqrt{\tau_{2}-\tau_1}\ind_{\{R_b\ge2\}}\right] =  \EL\left[R_b\sqrt{\tau_1}\sqrt{\tau_{2}-\tau_1}\ind_{\{2\le R_b< \gk \sqrt{L\log L}\}}\right]+o(1)} \nonumber\\ 
 & \hspace{1cm}\leq & \sqrt{L} (\log L)^{1/12} \EL\left[\sqrt{\tau_1}\sqrt{\tau_{2}-\tau_1}\ind_{\{R_b\ge2\}}\right] \nonumber\\
 & \hspace{1cm}& \indent + \gk \sqrt{L\log L} \EL\left[\sqrt{\tau_1}\sqrt{\tau_{2}-\tau_1}\ind_{\{R_b> \sqrt{L} (\log L)^{1/12}\}}\right] +o(1). \label{2eme partie} 
\end{eqnarray}
Let us deal with the second term:
\begin{eqnarray}
\lefteqn{\EL\left[\ind_{\{R_b> \sqrt{L} (\log L)^{1/12}\}}\sqrt{\tau_1}\sqrt{\tau_{2}-\tau_1}\right] }  \nonumber\\
  & = & \frac{1}{\P(b\in\tau)} \sum_{i=1}^b \sum_{j=1}^{b-i} \sqrt{i}\sqrt{j} \P\left(\tau_{1}=i,\tau_2-\tau_1=j,b\in\tau, R_b> \sqrt{L} (\log L)^{1/12}\right) \nonumber\\
  & = & \frac{1}{\P(b\in\tau)} \sum_{i=1}^b \sum_{j=1}^{b-i} \sqrt{i}\sqrt{j} K(i)K(j) \P\left(b-i-j\in\tau ,  
R_{b-i-j}> \sqrt{L} (\log L)^{1/12}-2\right). \label{second term part1}
\end{eqnarray}
But we have
\begin{eqnarray}
\lefteqn{\P\left( R_{b-i-j}> \sqrt{L} (\log L)^{1/12}  -2\left|b-i-j\in\tau \right.\right)  \leq  2 \P\left( R_{(b-i-j)/2}>  \frac{1}{2}\sqrt{L} (\log L)^{1/12}-1 \left|b-i-j\in\tau \right.\right) }\nonumber\\
 & \hspace{5cm}\leq & c \P\left( R_{(b-i-j)/2}>  \frac{1}{2}\sqrt{L} (\log L)^{1/12}-1 \right) \nonumber\\
 & \hspace{5cm}\leq & c \P\left(R_L>  \frac{1}{2}\sqrt{L} (\log L)^{1/12}-1\right) \leq c'\, e^{-c(\log L)^{1/6} } \label{second term part2}
\end{eqnarray}
where we first used Lemma~\ref{decond2} to take the conditioning off, and then Lemma~\ref{Wn grand}.
Putting~(\ref{second term part1}) and~(\ref{second term part2}) together, we get
\begin{eqnarray}
\lefteqn{ \EL\left[\ind_{\{R_b> \sqrt{L} (\log L)^{1/12}\}}\sqrt{\tau_1}\sqrt{\tau_{2}-\tau_1}\right]  }\nonumber\\
 & \hspace{2cm}\leq &  c' e^{-c(\log L)^{1/6} } \frac{1}{\P(b\in\tau)} \sum_{i=1}^b \sum_{j=1}^{b-i} \sqrt{i}\sqrt{j} K(i)K(j) \P\left(b-i-j\in\tau\right) \nonumber\\ 
 & \hspace{2cm}= & c' e^{-c(\log L)^{1/6} }\EL\left[\sqrt{\tau_1}\sqrt{\tau_{2}-\tau_1}\ind_{\{R_b\ge 2\}}\right].
\end{eqnarray}
So, recalling~(\ref{2eme partie}), we have
\begin{equation}
\EL\left[R_b\sqrt{\tau_1}\sqrt{\tau_{2}-\tau_1}\ind_{\{R_b\ge 2\}}\right] \leq 2\sqrt{L}(\log L)^{1/12}\EL\left[\sqrt{\tau_1}\sqrt{\tau_{2}-\tau_1}\ind_{\{R_b\ge2\}}\right]  +o(1)
\label{3eme partie}
\end{equation}
and we only have to estimate (recall \eqref{eq:doney})
\begin{eqnarray}
\EL\left[\sqrt{\tau_1}\sqrt{\tau_{2}-\tau_1}\ind_{\{R_b\ge2\}}\right] & =  &  \sum_{p=1}^{b-1} \sum_{q=1}^{b-p} \sqrt{p}\sqrt{q} K(p)K(q) \frac{\P(b-p-q\in\tau)}{\P(b\in\tau)} \nonumber\\
 & \leq & c \sqrt b \sum_{p=1}^{b-1} \sum_{q=1}^{b-p} \frac{1}{p\, q}\frac{1}{\sqrt{b+1-p-q}}. \label{Etau sqrt(gdt1)}
\end{eqnarray}
Using twice the elementary estimate 
\begin{eqnarray*}
\sum_{k=1}^{M-1} \frac{1}{k}\frac{1}{\sqrt{M-k}} \le c \frac{1}{\sqrt{M}} \log M,
\end{eqnarray*}
we get
\begin{eqnarray}
  \EL\left[\sqrt{\tau_1}\sqrt{\tau_{2}-\tau_1}\ind_{\{R_b\ge 2\}}\right]  \leq c \sqrt b \sum_{p=1}^{b-1} \frac{1}{p}\frac{1}{\sqrt{b-p+1}} \log (b-p+1) \leq c \sqrt b \frac{1}{\sqrt b} (\log L)^2.
\end{eqnarray}
Together with~(\ref{3eme partie}), this proves the desired estimate \eqref{112}.

\qed

\subsection{Dimension $d=4$ (a sketch)}
\label{rem:d4}

As we mentioned just after Theorem \ref{th:main}, it is possible to
adapt the change-of-measure argument to prove non-coincidence of
quenched and annealed critical points in dimension $d=4$ for the
general walks of Assumption \ref{hypwright}, while the method of
Birkner and Sun \cite{BS08} does not seem to adapt easily much beyond
the simple random walk case.  We will not give details, but for the
interested reader we hint at the ``right'' change of measure which
works in this case.

The ``change of measure function'' $g_{\cI}(Y)$ is still of the form
\eqref{eq:change1block}, factorized over the blocks which belong to 
$\mathcal I$, but this time  $M$ is a matrix with a finite
bandwidth:
\begin{eqnarray}
  F_k(Y)=-\frac c{\sqrt L}\sum_{i=L(k-1)+1}^{k L-p_0}\Delta_i\cdot
\Delta_{i+p_0},
\end{eqnarray}
where $p_0$ is an integer. The role of the normalization $L^{-1/2}$ is to
guarantee that Lemma \ref{borneM} still holds, and $p_0$ is to be chosen
such that $A(p_0)>0$, where $A(\cdot)$ is the function defined in Lemma
\ref{th:Ar}.
The existence of such $p_0$ is guaranteed by
the asymptotics \eqref{A(r)}, whose proof for $d=4$ is the same as for
$d=3$.

For the rest, the scheme of the proof of $\beta_c\ne\beta^{ ann}_c$
(in particular, the coarse-graining procedure) is analogous to that we
presented for $d=3$, and the computations involved are considerably
simpler.

\appendix
\section{Some technical estimates}

\begin{lemma} (Lemma A.2 in \cite{GLT08}) Let $\bP$ be the law of a
  recurrent renewal whose inter-arrival law satisfies
  $K(n)\stackrel{n\to\infty}\sim c_K n^{-3/2}$ for some $c_K>0$.
  There exists a constant $c>0$, that depends only on $K(\cdot)$, such
  that for any non-negative function $f_N(\tau)$ which depends only on
  $\tau\cap\{1,\ldots,N\}$, one has
\begin{equation}
\sup_{N>0} \frac{\E[f_N(\tau) \left|2N\in\tau\right.]}{\E[f_N(\tau)]}\leq c.
\end{equation}
\label{decond2}
\end{lemma}

\begin{lemma}
\label{th:app1} Under the same assumptions as in Lemma \ref{decond2}, and with
 $R_N:=|\tau\cap\{1,\ldots,N\}|$, there exists a constant $c>0$, such that for any positive
  function $\ga(N)$ which diverges at infinity and such as
  $\ga(N)=o(\sqrt N)$, we have
\begin{equation}
\P\left(R_N\geq \sqrt{N}\ga(N) \right) \leq e^{-c \ga(N)^2 + o\left(\ga(N)^2\right)} .
\end{equation}
\label{Wn grand}
\end{lemma}
{\sl Proof.}
For every $\gl>0$
\begin{eqnarray}
  \P\left(R_N\geq \sqrt{N}\ga(N) \right) & = & 
  \P\left(\tau_{\sqrt{N}\ga(N)}\leq N \right)  = 
  \P\left(\gl \ga(N)^2 \frac{\tau_{\sqrt{N}\ga(N)}}N\leq \gl \ga(N)^2 \right) \\
\nonumber
  & \leq & e^{\gl \ga(N)^2 } \E \left[e^{-\gl \frac{\ga(N)^2}N \tau_{\sqrt{N} \ga(N)}} \right]
  = e^{\gl \ga(N)^2 } \E \left[e^{-\gl \ga(N)^2\frac{\tau_1}{N} }\right]^{\sqrt{N}\ga(N)}. 
\end{eqnarray}
The asymptotic behavior of $\E \left[e^{-\gl \ga(N)^2\frac{\tau_1}{ N} }\right]$ is easily obtained:
\begin{eqnarray}
1-\E \left[e^{-\gl \ga(N)^2\frac{\tau_1}{ N } }\right] & = & \sum_{n\in\N} K(n)\left(1-e^{-n \gl \ga(N)^2 /N}\right)
\nonumber\\
&  \stackrel{N\to \infty}\sim&  c\frac{\sqrt{\gl} \ga(N)}{ \sqrt{N}} , \indent \indent \  c=c_K\int_0^{\infty} \frac{1-e^{-x}}{x^{3/2}} \dd x, \label{laplace transform tau} 
\end{eqnarray}
where 
the condition $\ga(N)^2/N \to 0$ was used to transform the sum into
an integral. Therefore, we get
\begin{eqnarray}
\E \left[e^{-\gl \ga(N)^2 \frac{\tau_1}{N} }\right]^{\sqrt{N}\ga(N)} & = &  \left(1-c\frac{\sqrt{\gl} \ga(N)}{ \sqrt{N}}+o\left(\frac{\ga(N)}{ \sqrt{N}}\right)\right)^{ \sqrt{N}\ga(N) } \nonumber\\
 & = & e^{- c \sqrt{\gl} \ga(N)^2 +o\left(\ga(N)^2\right)}.
\end{eqnarray}
Then, for any $\gl>0$,
\begin{equation}
\P\left(R_N\geq \sqrt{N}\ga(N) \right) \leq e^{(\gl- c \sqrt{\gl}) \ga(N)^2 +o\left(\ga(N)^2\right)}
\end{equation}
and taking $\gl=c^2/4$ we get the desired bound.
\qed

We need also the following standard result (cf. for instance \cite[Section 5]{GLT09}):
\begin{lemma}
Under the same hypothesis as in Lemma \ref{decond2}, we have the following convergence in law:
\begin{equation}
 \frac{c_K}{\sqrt{N}} R_{N} \stackrel{N\to\infty}{\Rightarrow} \frac{1}{\sqrt{2\pi}} \left| \cZ \right| \hspace{2cm} (\cZ\sim \cN(0,1)).
\end{equation}
\label{conv in law W}
\end{lemma}

\subsection{Proof of \eqref{eq:ijl}}
\label{app:ijl}

We wish to show that for distinct $i,j,l$ smaller than $r$, 
\begin{eqnarray}
\label{eq:A8}
  \bbE^X\bbE^Y[(\Delta_i\cdot\Delta_j)(\Delta_i\cdot\Delta_l)|X_r=Y_r]\le \frac cr.
\end{eqnarray}
We use the same method as in Section \ref{lemma EctF big}: we fix $x\in\Z^d$,
and we use the  notation $\Eyrx[\cdot] = \Ey[\cdot\left|Y_r=x\right.]$. Then,
\begin{eqnarray}
\Eyrx\left[(\gD_i\cdot\gD_j)(\gD_i\cdot\gD_l)\right] &  = &  \Eyrx\left[ (\gD_i\cdot\gD_j)\left(\gD_i\cdot\Eyrx\left[ \gD_l  \left|\gD_i,\gD_j\right.\right]\right)  \right] \nonumber\\
 & = & \frac{1}{r-2} \Eyrx\left[ (\gD_i\cdot\gD_j)\left(\gD_i\cdot(x-\gD_i-\gD_j)\right)\right] \nonumber\\
 & = & \frac{1}{r-2} \Eyrx\left[ (\gD_i\cdot\gD_j)\left((x\cdot\gD_i) -\norm{\gD_i}^2\right) - (\gD_i\cdot\gD_j)^2\right]\nonumber\\
 & \leq & \frac{1}{r-2} \Eyrx\left[ \left((x\cdot\gD_i) -\norm{\gD_i}^2\right)\left(\gD_i\cdot\Eyrx\left[ \gD_j  \left|\gD_i\right.\right]\right)  \right] \nonumber\\
 & =& \frac{1}{(r-1)(r-2)}  \Eyrx\left[ \left((x\cdot\gD_i) -\norm{\gD_i}^2\right)^2 \right] \nonumber \\
 & \leq &
\nonumber \frac{2}{(r-1)(r-2)} \Eyrx\left[\norm{x}^2\norm{\gD_i}^2 + 
\norm{\gD_i}^4\right]
\end{eqnarray}
and we can take by symmetry $i=1$.
Therefore,
\begin{eqnarray}
\label{Ect gDijk}
\lefteqn{\Ex\Ey\big[ \left.(\gD_i\cdot\gD_j)(\gD_i\cdot\gD_l)
\right|X_r=Y_r \big] =
\frac{   \Ex\big[\Ey\left[(\gD_i\cdot\gD_j)(\gD_i\cdot\gD_l)\left|Y_r=X_r\right.\right] \Py(Y_r=X_r) \big]}{ \Pxy(Y_r=X_r) 
} } \\\nonumber
&\hspace{2.5cm}\le & \frac{c}{r^2}\frac{\Ex\left[\left(
\|X_r\|^2\Ey\big[ \|\Delta_1\|^2|Y_r=X_r\big] + \Ey(\|\Delta_1\|^4|Y_r=X_r)
\right)
\Py(Y_r=X_r)
\right]}{\Pxy(Y_r=X_r)}.
\end{eqnarray}
At this point, as in the computations leading to \eqref{asympA part1}-\eqref{asympA
  part2}, one first notes that values $\|X_r\|\ge r^{3/5}$ or
$\|\Delta_1\|\ge r^{1/4}$ can be neglected; in the remaining range one
applies the Local Limit Theorem (both in the numerator and in the
denominator) and finally the sums are transformed
into  integrals. The estimate \eqref{eq:ijl} then follows after a few
lines of computation.  \qed

\subsection{Proof of \eqref{eq:ijkl}}
\label{app:ijkl}
We wish to prove that, for distinct $i,j,k,l\le r$,
\begin{equation}
\Ect[(\gD_i\cdot\gD_j)(\gD_k\cdot\gD_l)] \leq \frac{c}{r^2}.
\end{equation}
The proof is very similar to that of \eqref{eq:A8}, so we skip details. 
What one gets is that
\begin{equation}
 \Ect\left[ (\gD_i\cdot\gD_j)(\gD_k\cdot\gD_l) \right]   \leq 
\frac{c}{r^2}  \frac{ \Ex\left[\Ey\left[\left.Q\left(\frac{\norm{X_r}}{r^{1/2}},\{
\|\gD_i\|\}_{i=1,2,3}\right)\right| Y_r=X_r\right] \Py(Y_r=X_r)\right]}{
\Pxy(Y_r=X_r)},
\end{equation}
where $Q$ is a polynomial of degree $4$ in the variable $\|X_r\|/\sqrt r$.
Again, like after \eqref{Ect gDijk}, one uses the Local Limit Theorem to get the desired result.

\subsection{Proof of \eqref{ijklmn}}
\label{app:ijklmn}
In view of \eqref{dif intervals}, it suffices to prove that for $0<i\ne k\le r,0<j\ne l\le s$
\begin{eqnarray}
\label{rs}
\sum_{\nu,\mu=1}^3  \bbE^X\Ey[\Delta^{(\nu)}_i\Delta^{(\mu)}_k|X_r=Y_r]
\bbE^X\Ey[\Delta^{(\nu)}_j\Delta^{(\mu)}_l|X_s=Y_s]\le \frac c{rs}.
\end{eqnarray}
Both factors in the left-hand side have already been computed in  \eqref{454}-\eqref{455}.
Using these two expressions and once more the Local Limit Theorem, one arrives easily to \eqref{rs}.

\section*{Acknowledgments}
F. T. would like  to thank Rongfeng Sun for several enlightening discussions, 
and for showing us the preprint \cite{BS09}
before its publication.

\end{document}